\RequirePackage{fix-cm}
\documentclass[envcountsame,smallextended]{svjour3}       
\smartqed  
\usepackage{amsthm}
\usepackage{amssymb}
\usepackage{graphicx}
\usepackage[top=35mm, bottom=25mm, left=18mm, right=18mm]{geometry}
\usepackage[numbers,sort&compress]{natbib}

\usepackage{booktabs}
\usepackage{multirow}
\usepackage{orcidlink}
\usepackage{tabularx}

\usepackage{amsmath}
\usepackage{amsfonts}
\usepackage{amssymb}
\usepackage{enumitem}
\usepackage{graphicx}
\usepackage{subcaption}
\usepackage{tikz}
\usetikzlibrary{arrows.meta, positioning}
\usepackage{pgfplots}
\usetikzlibrary{patterns}
\usepackage{multicol}
\usepackage{hyperref}
\usepackage{xcolor}
\usepackage{cancel}
\usepackage{ulem}
\usepackage[active]{srcltx}

\numberwithin{equation}{section}
\usepackage{subcaption}

\usepackage{booktabs}
\usepackage{multirow}

\usepackage{algorithm}%
\usepackage{algorithmicx}%
\usepackage{algpseudocode}%

\newcommand{\MRAR}{\mathrm{MR}\text{-}\mathrm{AR}}
\newcommand{\sym}{\operatorname{Sym}}
\newcommand{\cF}{\mathcal{F}}
\newcommand\R{\mathbb{R}}
\newcommand{\N}{\mathbb{N}}

\newcommand{\eps}{\epsilon}

\newcommand\ov[1]{\overline{#1}}

\newcommand\bs[1]{\boldsymbol{#1}}

\newcommand\Nz{\mathbb{N}_0}

\newcommand{\norm}[1]{\left\lVert#1\right\rVert}

\newcommand\inft[1]{\mathop{\bs{\inf}}\limits_{#1}}

\newcommand\maxt[1]{\mathop{\bs{\max}}\limits_{#1}}
\newcommand\mint[1]{\mathop{\bs{\min}}\limits_{#1}}
\let\oldliminf\liminf
\renewcommand{\liminf}{\bs{\oldliminf}}
\let\oldlimsup\limsup
\renewcommand{\limsup}{\bs{\oldlimsup}}
\let\oldmax\max
\renewcommand{\max}{\bs{\oldmax}}
\let\oldmin\min
\renewcommand{\min}{\bs{\oldmin}}
\let\oldsup\sup
\renewcommand{\sup}{\bs{\oldsup}}
\let\oldinf\inf
\renewcommand{\inf}{\bs{\oldinf}}
\let\oldlim\lim
\renewcommand{\lim}{\bs{\oldlim}}

\journalname{}

\begin{document}

\title{Multi-rate high-order regularization with an order-$(p-2)$ oracle
}

\titlerunning{Multi-rate high-order regularization}        

\author{Alireza Kabgani\,\orcidlink{0000-0001-6554-6969}
}

\authorrunning{A. Kabgani} 

\institute{A. Kabgani\at
              Department of Mathematics, University of Antwerp, Antwerp, Belgium. \\
              \email{alireza.kabgani@uantwerp.be}        
}

\date{Received: date / Accepted: date}

\maketitle

\vspace{-4mm}
\begin{abstract}
We consider unconstrained minimization of a possibly nonconvex function $f:\R^n\to\R$ whose $p$th derivative is Lipschitz continuous, with $p\geq3$, when derivative information is available only through order $p-2$. The unavailable tensors $\nabla^{p-1}f$ and $\nabla^pf$ are reconstructed from central finite differences of $\nabla^{p-2}f$. Since the two reconstructions have different point-query costs, we use a block length $m$ for the order-$p$ approximation and a subblock length $\ell\leq m$ for the order-$(p-1)$ approximation. The latter is transported between consecutive reconstruction points using the stored order-$p$ approximation. This leads to an adaptive regularization method that does not require prior knowledge of the Lipschitz constant. We prove finite termination at an $\eps$-approximate first-order stationary point and derive an explicit bound on the number of bundled order-$(p-2)$ oracle queries. For $m=(p-2)n^2$ and $\ell=(p-2)n$, the accuracy-dependent term has dimension--accuracy scaling $n^{3/p}\eps^{-(p+1)/p}$. Thus the method retains the $\eps^{-(p+1)/p}$ dependence of full order-$p$ regularization while using derivatives only through order $p-2$. Within the dense reconstruction framework considered here, using separate block and subblock lengths gives dimension dependence of order $n^{3/p}$, compared with $n^{4/p}$ when a common length is imposed.


\end{abstract}

\section{Introduction}
\label{intro}

We consider the unconstrained minimization problem
\begin{equation}\label{eq:problem}
\mint{x\in\R^n} f(x),
\end{equation}
where $f:\R^n\to\R$ may be nonconvex. Throughout the paper, $p\geq 3$ is an integer, $f\in C^p(\R^n)$ is bounded below, and $\nabla^pf$ is globally Lipschitz continuous with constant $L_p>0$. The precise assumptions and notation are given in Section~\ref{sec:preliminaries}.
Given $\eps>0$, our goal is to compute a point $x\in\R^n$ satisfying $\norm{\nabla f(x)}\leq\eps$. We refer to such a point as an $\eps$-approximate first-order stationary point \cite{Nesterov2018}.

Regularized Taylor methods exploit higher-order smoothness to improve the worst-case dependence on the target accuracy when the corresponding derivative information is available. Cubic regularization achieves an evaluation-complexity bound of order $\mathcal O(\eps^{-3/2})$ under Lipschitz continuity of the Hessian \cite{CartisGouldToint2011ARC,griewank1981modification,NesterovPolyak2006}. 
More generally, when $\nabla^pf$ is Lipschitz continuous, regularized Taylor methods of degree $p$ attain the first-order evaluation complexity
$\mathcal O(\eps^{-(p+1)/p})$ \cite{BirginEtAl2017}. Related high-order Taylor models and their associated tensor-envelope constructions in nonconvex composite
optimization are studied in \cite{Kabgani2026Tensor}.
 Sharp evaluation-complexity results for higher-order optimality measures were established in \cite{CartisGouldToint2020Sharp}, while the $\eps^{-(p+1)/p}$ dependence for first-order stationarity is also a lower bound in an oracle model that permits access to derivatives of all orders \cite{CarmonEtAl2020}.
 
 These methods require derivative information through order $p$. In many problems, however, the highest-order derivatives may be unavailable or expensive to evaluate. A natural question is whether, when an oracle provides derivatives only through order $q<p$, one can still exploit the $p$th-order smoothness of the objective and retain the $\eps^{-(p+1)/p}$ dependence on the target accuracy. This question is particularly challenging when more than one of the highest-order Taylor tensors is unavailable.

Finite differences provide a natural way to approximate unavailable derivatives from lower-order information. For cubic regularization, finite-difference Hessian approximations have been used to obtain first-order implementations with oracle complexity $\mathcal O(n\eps^{-3/2})$ \cite{GrapigliaGoncalvesSilva2022}. Finite-difference approximations have also been studied for quadratic regularization \cite{Grapiglia2023}, while adaptive regularization frameworks allowing inexact function and derivative evaluations were developed in \cite{BellaviaEtAl2019}. Another way to reduce derivative cost is to reuse previously computed approximations over several iterations. Lazy exact-Hessian updates were studied in \cite{DoikovChaytiJaggi2023}, and combining finite-difference approximations with lazy Hessian reuse yields first-order and zeroth-order oracle complexities $\mathcal O(n^{1/2}\eps^{-3/2})$ and $\mathcal O(n^{3/2}\eps^{-3/2})$, respectively \cite{DoikovGrapiglia2025}. Related adaptive reuse with inexact gradients is considered in \cite{GoncalvesGehlen2026}. At higher orders, secant updates for derivative tensors were developed in \cite{WelzelHauser2024}, while \cite{CartisJerad2025} studies adaptive regularization with lazy, finite-difference, and secant approximations of the highest-order tensor.

The case $q=p-1$, where only the highest-order Taylor tensor is unavailable, has been studied recently in
\cite{CartisJerad2025,DoikovGrapiglia2026}. In particular, \cite{DoikovGrapiglia2026} reconstructs $\nabla^pf$ from finite differences of $\nabla^{p-1}f$ and reuses the approximation over blocks of iterations. With an appropriate block length, the resulting method requires $\mathcal O_p\!\left(n^{1/p}\eps^{-(p+1)/p}\right)$ order-$(p-1)$ oracle calls in the accuracy-dependent term. The same work leaves the case $q<p-1$, for $p\geq3$, as an open direction and points out the difficulty of constructing multivariate finite-difference approximations with explicit error bounds suitable for controlling the accumulated approximation error.

We consider the case $q=p-2$. Thus, an oracle query at $x$ provides $\{f(x),\nabla f(x),\ldots,\nabla^{p-2}f(x)\}$, 
while the two consecutive Taylor tensors $\nabla^{p-1}f(x)$ and $\nabla^p f(x)$ are unavailable. In contrast
with the case $q=p-1$, both tensors must now be reconstructed from lower-order information. Moreover, when the approximations are reused away from their reconstruction points, the error in the order-$p$ approximation also enters the transport of the order-$(p-1)$ approximation.

The two reconstructions have different point-query costs. The constructions in Section~\ref{sec:reconstruction} require $2n$ shifted bundled order-$(p-2)$ oracle queries to reconstruct $\nabla^{p-1}f$, but $n(n+1)$ shifted queries to reconstruct $\nabla^pf$. Reconstructing both tensors at every model step would therefore incur a quadratic dimension dependence per step. This asymmetry motivates reconstructing the two tensors on different schedules.

We therefore use a block length $m$ and a subblock length $\ell\leq m$. An approximation of $\nabla^pf$ is reconstructed at the beginning of each block and reused throughout it, while the approximation of $\nabla^{p-1}f$ is reconstructed at the beginning of each subblock. Between two consecutive subblock reconstructions,
the latter is transported using the stored order-$p$ approximation. The resulting transport error depends on both the displacement from the block reference and the displacement from the current subblock reference, leading to two distinct path scales in the analysis.

Optimizing the two lengths in the dense point-query bound gives the dimension factor $\Theta_p(n^{3/p})$ in the accuracy-dependent term, whereas imposing a common length gives $\Theta_p(n^{4/p})$.

\subsection{Contributions}
\label{sec:contributions}

The main contributions of this work are as follows.

\begin{enumerate}[label=\textup{(\roman*)},leftmargin=*]

\item \textbf{Reconstruction of the two missing Taylor tensors.}
We construct finite-difference approximations of both $\nabla^{p-1}f$ and $\nabla^pf$ from evaluations of $\nabla^{p-2}f$. The first reconstruction uses central differences, while the second combines central second differences, polarization, coordinate assembly, and symmetrization. We establish explicit operator-norm error bounds and exact bundled-oracle query counts for both constructions.

\item \textbf{Multi-rate adaptive regularization.}
We develop $\MRAR(p,p-2)$ using a block length $m$ for the order-$p$ approximation and a subblock length $\ell\leq m$ for the order-$(p-1)$ approximation. Between consecutive subblock reconstructions, the latter is transported using the stored order-$p$ approximation. We establish the corresponding two-reference transport bound and control the resulting block and subblock displacement terms after summation over a complete block. The method adapts to the unknown Lipschitz constant $L_p$.

\item \textbf{Bundled-oracle complexity.}
We prove finite termination at an $\eps$-approximate first-order stationary point and give an explicit bound on the number of bundled order-$(p-2)$ oracle queries. For $m=(p-2)n^2$ and $\ell=(p-2)n$, the method computes a point satisfying $\|\nabla f(x)\|\leq\eps$ using
\begin{equation}\label{eq:intro-main-bound}
\mathcal{O}_p\!\left(n^{3/p}L_{\max}^{1/p}\Delta_f\eps^{-(p+1)/p}+n^2\left[1+\log_2\!\left(\frac{L_{\max}}{L_0}\right)\right]\right),
\end{equation}
bundled-oracle point queries, where $\Delta_f=f(x_0)-f_{\rm low}$ and $L_{\max}=\max\{L_0,2L_p\}$. Thus the accuracy-dependent term retains the $\eps^{-(p+1)/p}$ dependence of full order-$p$ regularization while using derivatives only through order $p-2$. Within the dense reconstruction framework considered here, optimizing the two lengths gives the dimension factor $\Theta_p(n^{3/p})$, compared with $\Theta_p(n^{4/p})$ when a common length is imposed.
\end{enumerate}
Here and throughout, the constants implicit in $\mathcal O_p(\cdot)$ and $\Theta_p(\cdot)$ depend only on $p$. The complexity measure in
\eqref{eq:intro-main-bound} counts point queries to the bundled oracle. Arithmetic operations, storage of dense tensors, and the work required to compute an admissible model step are not included in this oracle count and are discussed separately in Section~\ref{sec:complexity}.

\subsection{Organization}
The remainder of the paper is organized as follows.
Section~\ref{sec:preliminaries} introduces the assumptions, notation, bundled-oracle model, and Taylor estimates. Section~\ref{sec:reconstruction} develops the two finite-difference reconstructions and the affine transport bound, and Section~\ref{sec:algorithm} uses these ingredients to define the multi-rate adaptive regularization method. Section~\ref{sec:analysis} establishes the blockwise error control, sufficient decrease, and finite termination. Section~\ref{sec:complexity} derives the bundled order-$(p-2)$ oracle complexity, studies the choice of the block and subblock lengths, and discusses computational costs not captured by the point-query measure.
Section~\ref{sec:conclusion} concludes the paper. The appendices contain the choice of constants, auxiliary estimates, and technical proofs.



\section{Preliminaries and notation}
\label{sec:preliminaries}

The space $\R^n$ is equipped with the Euclidean norm $\Vert\cdot\Vert$ and inner product $\langle\cdot,\cdot\rangle$. Let $e_1,\ldots,e_n$ denote the standard orthonormal basis of $\R^n$, and let $(u)_i$ denote the $i$th component of $u\in\R^n$. We write $\Nz:=\N\cup\{0\}$. 

For nonnegative quantities $A$ and $B$, we write $A=\mathcal O_p(B)$ if there exists a constant $C_p>0$, depending only on $p$, such that $A\leq C_pB$ uniformly in all other quantities appearing in the bound. Likewise, $A=\Theta_p(B)$ means that there exist constants $0<c_p\leq C_p<\infty$, depending only on $p$, such that $c_pB\leq A\leq C_pB$ uniformly in all other quantities.

 For an integer $j\geq1$ and a $j$-linear form $T:(\R^n)^j\to\R$, not necessarily symmetric, define
\[
\Vert T\Vert_{[j]}:=\maxt{\substack{\Vert u_i\Vert\leq1\\ i=1,\ldots,j}}\vert T[u_1,\ldots,u_j]\vert.
\]
For $s\in\R^n$, we write
\[
T[s]^j:=T[\underbrace{s, \ldots, s}_{j\text{ times}}].
\]
More generally, for $1\leq r<j$, $T[u_1, \ldots, u_r]$ denotes the $(j-r)$-linear form defined by
\[
\left(T[u_1, \ldots, u_r]\right)[v_1, \ldots, v_{j-r}]:=T[u_1, \ldots, u_r, v_1, \ldots, v_{j-r}].
\]
For a $j$-linear form $T$, define its symmetrization by
\[
\sym(T) := \frac{1}{j!}\sum_{\pi\in \mathcal{S}_j} T^\pi,
\]
where $\mathcal{S}_j$ is the permutation group of $\{1,\ldots,j\}$ and
\[
T^\pi[u_1, \ldots, u_j] := T[u_{\pi(1)}, \ldots, u_{\pi(j)}].
\]
Then $\sym(T)$ is symmetric and satisfies $\sym(T)[s]^j=T[s]^j$ and $\Vert \sym(T)\Vert_{[j]}\leq \Vert T\Vert_{[j]}$.
For a $k$-linear form $E$ and $e\in \R^n$, define the $(k+1)$-linear form $E\otimes e$ by
\[
\left( E\otimes e\right)\left[u_1, \ldots, u_{k+1}\right] := E\left[u_1, \ldots, u_{k}\right]  \langle e, u_{k+1}\rangle.
\]

\begin{assumption}[Standing assumptions]\label{ass:main}
Let $p\geq3$ be an integer and $f\in C^p(\R^n)$. Assume that there exist $f_{\rm low}\in\R$ and $L_p>0$ such that $f(x)\geq f_{\rm low}$ for all $x\in\R^n$,  and
\begin{equation}\label{eq:Lp}
\Vert\nabla^p f(x)-\nabla^p f(y)\Vert_{[p]} \leq L_p\Vert x-y\Vert, \qquad \forall x,y\in\R^n.
\end{equation}
\end{assumption}

\begin{definition}[Bundled order-$(p-2)$ oracle]\label{def:oracle}
For $x\in\R^n$, define the bundled order-$(p-2)$ oracle by
\[
\cF(x) := \left(f(x),\nabla f(x),\ldots,\nabla^{p-2}f(x) \right).
\]
A bundled order-$(p-2)$ oracle query at $x$ returns the exact bundle $\mathcal{F}(x)$ and is counted as one point query.
\end{definition}
Previously retained bundles are reused without an additional oracle query. Any evaluation separately prescribed by the algorithm or by a finite-difference reconstruction is counted as a new query, even if its query point coincides with a previously evaluated point.

For a target accuracy $\eps>0$, let $N_{p-2}(\eps)$ denote the total number of bundled order-$(p-2)$ point queries made up to and including termination. This count does not include arithmetic operations, dense tensor storage and assembly, or the work required by the inner model solver. These costs are discussed in Section~\ref{sec:complexity}.

Since $f\in C^p(\R^n)$, the derivatives $\nabla^jf(x)$, $1\leq j\leq p$, are symmetric $j$-linear forms. In particular, 
\[
D\nabla^{p-2}f(x)[v]=\nabla^{p-1}f(x)[v],
\qquad
D^2\nabla^{p-2}f(x)[u,v]=\nabla^p f(x)[u,v].
\]
 For $x,s\in\R^n$, define the degree-$p$ Taylor polynomial of $f$ at $x$ by 
\begin{equation}\label{eq:Taylor}
T_p(x,s) := f(x) + \sum_{j=1}^{p} \frac{1}{j!}\nabla^jf(x)[s]^j. 
\end{equation}
Its gradient with respect to the displacement $s$ is
\begin{equation}\label{eq:Taylor-gradient}
\nabla_sT_p(x,s) = \nabla f(x) + \sum_{j=2}^{p} \frac{1}{(j-1)!} \nabla^jf(x)[s]^{j-1}.
\end{equation}
Taylor's theorem and \eqref{eq:Lp} give, for all $x,s\in\R^n$,
\begin{equation}\label{eq:function-rem}
\left|f(x+s)-T_p(x,s)\right|\leq \frac{L_p}{(p+1)!}\Vert s\Vert^{p+1},
\end{equation}
and
\begin{equation} \label{eq:gradient-rem}
\Vert\nabla f(x+s)-\nabla_sT_p(x,s)\Vert\leq \frac{L_p}{p!}\Vert s\Vert^{p}.
\end{equation}

\section{Reconstruction and affine transport}
\label{sec:reconstruction}
This section develops the finite-difference reconstructions of the two unavailable tensors and the affine transport estimate used between reconstruction points.

Applying Taylor's formula with integral remainder to $x\mapsto \nabla^{p-2} f(x)$ gives, for $x, v\in \R^n$,
\begin{equation}\label{eq:lem:sec:expan:1}
\nabla^{p-2} f(x+v) =  \nabla^{p-2} f(x) + \nabla^{p-1} f(x)[v] +\frac{1}{2}\nabla^{p} f(x)[v,v] + R_x(v),
\end{equation}
where 
\begin{equation}\label{eq:lem:sec:expan:2}
R_x(v) :=  \int_0^1(1-t)\left[\nabla^{p} f(x+tv)- \nabla^{p} f(x)\right][v,v] dt.
\end{equation}
By \eqref{eq:Lp},
\begin{equation}\label{eq:lem:sec:expan:3}
\Vert R_x(v) \Vert_{[p-2]}\leq L_p \Vert v\Vert^3 \int_{0}^{1}t(1-t)dt = \frac{L_p}{6}\Vert v\Vert^3.
\end{equation}

\subsection{Finite-difference reconstructions}
We first reconstruct $\nabla^{p-1}f$. For $x, v\in \R^n$ and $h>0$, define
\begin{equation}\label{eq:DA}
D_A(x; v, h) :=\frac{\nabla^{p-2}f(x+hv) - \nabla^{p-2}f(x-hv)}{2h}.
\end{equation}
This is a symmetric $(p-2)$-linear form approximating $\nabla^{p-1}f(x)[v]$. Define
\[
\widetilde{A}_h(x) := \sum_{i=1}^{n} D_A(x; e_i, h) \otimes e_i.
\]
This construction extends the central-difference estimate of \cite[Lemma~5]{DoikovGrapiglia2025} to tensor-valued derivatives.
\begin{lemma}[Reconstruction of $\nabla^{p-1}f$]\label{lem:rec:p-1}
For every $x,v\in\R^n$ and $h>0$,
\begin{equation}\label{eq:lem:rec:p-1:1}
\Vert D_A(x; v, h) - \nabla^{p-1}f(x)[v]\Vert_{[p-2]}\leq \frac{L_p}{6}h^2 \Vert v\Vert^3,
\end{equation}
and
\begin{equation}\label{eq:lem:rec:p-1:2}
\Vert \sym(\widetilde{A}_h(x)) -  \nabla^{p-1}f(x)\Vert_{[p-1]} \leq \frac{\sqrt{n}L_p}{6} h^2.
\end{equation}
The reconstruction requires exactly $2n$ bundled order-$(p-2)$ oracle queries, at the shifted points $x\pm he_i$, $i=1,\ldots,n$.
\end{lemma}
\begin{proof}
Subtracting \eqref{eq:lem:sec:expan:1} with increments $hv$ and $-hv$ gives
\[
D_A(x; v, h) - \nabla^{p-1}f(x)[v] = \frac{R_x(hv)- R_x(-hv)}{2h},
\]
and hence \eqref{eq:lem:rec:p-1:1} follows from \eqref{eq:lem:sec:expan:3}. Since $\nabla^{p-1}f(x)$ is symmetric and multilinear, we have 
\[
\nabla^{p-1}f(x) = \sum_{i=1}^{n}\nabla^{p-1} f(x) [e_i]\otimes e_i.
\]
Let $E_i:=D_A(x; e_i, h) - \nabla^{p-1} f(x) [e_i]$. For unit vectors $u_1, \ldots, u_{p-1}\in \R^n$, we have
\begin{align*}
\left\vert \left(\widetilde{A}_h(x) - \nabla^{p-1}f(x)\right)[u_1, \ldots, u_{p-1}]\right\vert
&= \left\vert \sum_{i=1}^{n}\left(E_i\otimes e_i\right) [u_1, \ldots, u_{p-1}]\right\vert 
= \left\vert \sum_{i=1}^{n}E_i[u_1, \ldots, u_{p-2}](u_{p-1})_i \right\vert 
\\&\overset{(i)}{\leq} \left( \sum_{i=1}^{n}\left\vert E_i[u_1, \ldots, u_{p-2}] \right\vert^2\right)^{1/2}
\left(\sum_{i=1}^{n}\left\vert(u_{p-1})_i \right\vert^2\right)^{1/2}
\\&\overset{(ii)}{\leq} \left(\sum_{i=1}^{n} \Vert E_i\Vert_{[p-2]}^2\right)^{1/2}.
\end{align*}
where $(i)$ follows from Cauchy--Schwarz and  $(ii)$ uses $\Vert u_{p-1}\Vert=1$.
Taking the supremum and using \eqref{eq:lem:rec:p-1:1},  we obtain
\begin{align*}
\Vert \sym(\widetilde{A}_h(x)) -  \nabla^{p-1}f(x)\Vert_{[p-1]}&= \Vert \sym(\widetilde{A}_h(x) -  \nabla^{p-1}f(x))\Vert_{[p-1]}
\\&\leq \Vert \widetilde{A}_h(x) -  \nabla^{p-1}f(x)\Vert_{[p-1]}\leq \left(\sum_{i=1}^{n} \Vert E_i\Vert_{[p-2]}^2\right)^{1/2} \leq \frac{\sqrt{n}L_p}{6} h^2,
\end{align*}
which proves \eqref{eq:lem:rec:p-1:2}. Finally, each of the $n$ central differences in \eqref{eq:DA} uses the two shifted points $x\pm he_i$. Hence the reconstruction requires exactly $2n$ shifted oracle queries.
\end{proof}

We next reconstruct $\nabla^pf$. For $x,v\in\R^n$ and $h>0$, define
\begin{equation}\label{eq:D_B}
 D_B(x;v,h) := \frac{\nabla^{p-2} f(x+hv)-2\nabla^{p-2} f(x)+\nabla^{p-2} f(x-hv)}{h^2}.
\end{equation}
This is a symmetric $(p-2)$-linear form approximating $\nabla^pf(x)[v,v]$. 
For $p=3$, related central second differences of gradients are used in \cite{Nesterov2021Superfast} to approximate third-order directional information.

For $1\leq i<j\leq n$, set $u_{ij} := \frac{e_i+e_j}{\sqrt{2}}$. By symmetry of $\nabla^p f(x)$, the polarization identity \cite[Eq.~(1b)]{Thomas2014Polarization} yields
\begin{equation}\label{eq:polarization}
\nabla^p f(x)[e_i, e_j] = \nabla^p f(x)[u_{ij}, u_{ij}] - \frac{1}{2}\nabla^p f(x)[e_i, e_i] - \frac{1}{2}\nabla^p f(x)[e_j, e_j].
\end{equation}
Accordingly, set $Q_{ii} :=  D_B(x;e_i,h)$, $i=1, \ldots, n$, and, for $1\leq i< j\leq n$,
\begin{equation}\label{eq:Qij}
Q_{ij}=Q_{ji} := D_B(x;u_{ij},h) - \frac{1}{2}Q_{ii}-  \frac{1}{2}Q_{jj}.
\end{equation}
Thus, $Q_{ij}$ is a symmetric $(p-2)$-linear approximation of $\nabla^p f(x)[e_i, e_j]$. Define 
\begin{equation}\label{eq:B-raw}
\widetilde B_h(x)[v_1,\ldots,v_p] :=\sum_{i,j=1}^n (v_1)_i(v_2)_j Q_{ij}[v_3,\ldots,v_p].
\end{equation}

\begin{lemma}[Reconstruction of $\nabla^p f$]\label{lem:nab-p-reconstruction}
For every $x,v\in\R^n$ and $h>0$, 
\begin{equation}\label{eq:B-directional}
\norm{D_B(x;v,h)-\nabla^p f(x)[v,v]}_{[p-2]} \leq \frac{L_p}{3}h\norm{v}^3.
\end{equation}
Moreover,
\begin{equation}\label{eq:B-full}
\norm{\sym\left(\widetilde B_h(x)\right)-\nabla^p f(x)}_{[p]} \leq  \frac{2n}{3}L_ph.
\end{equation}
If $\cF(x)$ is already available, the reconstruction  requires exactly $n(n+1)$ additional bundled order-$(p-2)$ oracle queries.
\end{lemma}
\begin{proof}
Adding \eqref{eq:lem:sec:expan:1} with increments $hv$ and $-hv$ gives
\[
D_B(x;v,h) - \nabla^p f(x)[v,v] = \frac{R_x(hv)+R_x(-hv)}{h^2}.
\]
Hence, by \eqref{eq:lem:sec:expan:3},
\[
\Vert D_B(x;v,h) - \nabla^p f(x)[v,v]\Vert_{[p-2]} \leq \frac{L_p}{3} h\Vert v\Vert^3,
\]
which proves \eqref{eq:B-directional}. For $1\leq i<j\leq n$, equations \eqref{eq:polarization} and \eqref{eq:Qij} give
\begin{align*}
Q_{ij} - \nabla^p f(x)[e_i, e_j]& = D_B(x;u_{ij},h)-\frac{1}{2}Q_{ii}-\frac{1}{2}Q_{jj} -
\nabla^p f(x)[u_{ij}, u_{ij}] + \frac{1}{2}\nabla^p f(x)[e_i, e_i] + \frac{1}{2}\nabla^p f(x)[e_j, e_j]
\\
& = D_B(x;u_{ij},h)-\nabla^p f(x)[u_{ij}, u_{ij}]-\frac{1}{2}\left(Q_{ii}-\nabla^p f(x)[e_i, e_i]\right)-\frac{1}{2}\left(Q_{jj}-\nabla^p f(x)[e_j, e_j]\right).
\end{align*}
Since $\norm{e_i}=\norm{u_{ij}}=1$, \eqref{eq:B-directional} gives
\begin{align*}
\left\Vert Q_{ij} - \nabla^p f(x)[e_i, e_j] \right\Vert_{[p-2]} 
\leq \frac{L_p}{3} h +\frac{L_p}{6} h+\frac{L_p}{6} h = \frac{2L_p}{3} h.
\end{align*}
For $i=j$, \eqref{eq:B-directional} gives 
\[
\norm{Q_{ii}-\nabla^pf(x)[e_i,e_i]}_{[p-2]}\leq\frac{L_p}{3}h.
\]
Thus,
\[
\norm{Q_{ij}-\nabla^pf(x)[e_i,e_j]}_{[p-2]}\leq\frac{2L_p}{3}h\qquad\forall i,j.
\]
Now let $a,b,v_3,\ldots,v_p$ be unit vectors. By multilinearity,
\[
\nabla^p f(x)[a,b,v_3,\ldots,v_p] = \sum_{i,j=1}^n (a)_i (b)_j \nabla^p f(x)[e_i, e_j][v_3,\ldots,v_p].
\]
Consequently,
\begin{align*}
\left|(\widetilde B_h(x)-\nabla^p f(x))[a,b,v_3,\ldots,v_p]\right|&\leq \sum_{i,j=1}^n |(a)_i||(b)_j|\left\Vert Q_{ij} - \nabla^p f(x)[e_i, e_j] \right\Vert_{[p-2]}
\\&\leq \frac{2L_p}{3} h \sum_{i,j=1}^n |(a)_i||(b)_j| = \frac{2L_p}{3} h\norm{a}_1\norm{b}_1 \le n\frac{2L_p}{3} h,
\end{align*}
where the last inequality follows from $\norm{a}_1\le\sqrt n\norm{a}=\sqrt n$ and similarly for $b$.
Taking the supremum, this gives
\begin{align*}
\Vert \sym(\widetilde B_h(x)) -  \nabla^{p}f(x)\Vert_{[p]}= \Vert \sym(\widetilde B_h(x) -  \nabla^{p}f(x))\Vert_{[p]}
\leq \Vert \widetilde B_h(x) -  \nabla^{p}f(x)\Vert_{[p]}\leq \frac{2n}{3}L_ph,
\end{align*}
which proves \eqref{eq:B-full}. Finally, the reconstruction uses the $n$ coordinate directions $e_i$, $i=1,\ldots,n$, and one direction $u_{ij}$ for each pair
$1\leq i<j\leq n$. Hence the total number of directions is
\[
n+\binom{n}{2} = \frac{n(n+1)}{2}.
\]
For each direction $v$, \eqref{eq:D_B} requires evaluations at the two shifted points $x\pm hv$. Since $\nabla^{p-2} f(x)$ is already contained in $\cF(x)$, no additional evaluation at $x$ is needed. Hence the reconstruction requires $n(n+1)$
additional bundled-$(p-2)$ oracle queries.
\end{proof}

\subsection{Affine transport}
The two reconstructions above have different point-query costs. In the algorithm developed in the next section, we will therefore avoid reconstructing both tensors whenever the current point changes. We first show how an approximation of $\nabla^{p-1}f$ at $y$ can be transported to a point $x$ using an approximation of $\nabla^pf$ available at a possibly different point $z$.

Let $z\in\R^n$ be a point at which a symmetric $p$-linear approximation $\widehat B_z$ of $\nabla^p f(z)$ is available, and let $y\in\R^n$ be a point at which a symmetric $(p-1)$-linear approximation $\widehat A_y$ of $\nabla^{p-1} f(y)$ is available. For a current point $x\in\R^n$, define the affine transport
\begin{equation}\label{eq:transport-def}
\widehat A_x :=\widehat A_y+\widehat B_z[x-y].
\end{equation}
Then $\widehat A_x$ is a symmetric $(p-1)$-linear form.

\begin{lemma}[Affine transport bound]\label{lem:transport}
Suppose that
\[
\norm{\widehat A_y-\nabla^{p-1}f(y)}_{[p-1]} \leq \delta_{A}, \qquad \norm{\widehat B_z-\nabla^{p}f(z)}_{[p]} \leq \delta_{B},
\]
for some $\delta_{A},\delta_{B}\geq 0$. Then, for every $x\in\R^n$,
\begin{equation}\label{eq:transport}
\norm{\widehat A_x-\nabla^{p-1}f(x)}_{[p-1]} \leq \delta_{A} + \left(\delta_{B}+L_p\norm{x-z}\right)\norm{x-y} + \frac{L_p}{2}\norm{x-y}^2.
\end{equation}
\end{lemma}
\begin{proof}
Set $d:=x-y$. By the fundamental theorem of calculus,
\[
\nabla^{p-1}f(x)-\nabla^{p-1}f(y) = \int_0^1 \nabla^{p}f(y+td)[d]\,dt.
\]
Using $\widehat A_x  = \widehat A_y+\widehat B_z[d]$ and adding and subtracting $\nabla^{p}f(z)[d]$ gives
\[
\widehat A_x-\nabla^{p-1}f(x) = \widehat A_y-\nabla^{p-1}f(y) + (\widehat B_z-\nabla^{p}f(z))[d] + \int_0^1 \left(\nabla^{p}f(z)-\nabla^{p}f(y+td)\right)[d]\,dt.
\]
It follows that
\[
\norm{\widehat A_x-\nabla^{p-1}f(x)}_{[p-1]}\leq \delta_{A}+\delta_{B}\norm d + L_p\norm d \int_0^1\norm{z-(y+td)}\,dt.
\]
Note that $\norm{z-(y+td)} \leq (1-t)\norm{z-y}+t\norm{z-x}$.
Integrating over $t\in[0,1]$ yields
\[
\int_0^1 \norm{z-(y+td)}\,dt \leq \frac{1}{2} \left(\norm{y-z}+\norm{x-z}\right).
\]
Therefore,
\[
\norm{\widehat A_x-\nabla^{p-1}f(x)}_{[p-1]} \leq \delta_{A} + \delta_{B}\norm{x-y} + \frac{L_p}{2} \left(\norm{y-z}+\norm{x-z}\right)\norm{x-y}.
\]
Finally, $\norm{y-z} \leq \norm{x-y}+\norm{x-z}$ implies 
\begin{align*}
\norm{\widehat A_x-\nabla^{p-1}f(x)}_{[p-1]} & \leq \delta_{A} + \delta_{B}\norm{x-y} + \frac{L_p}{2} \left(\norm{x-y}+2\norm{x-z}\right)\norm{x-y},
\end{align*}
which proves \eqref{eq:transport}.
\end{proof}

The bound \eqref{eq:transport} contains two displacement scales. The factor $\Vert x-y\Vert$ measures movement since the most recent order-$(p-1)$ reconstruction, while $\Vert x-z\Vert$ measures movement from the point at which the stored order-$p$ tensor was reconstructed. This interaction is the source of the two-scale path terms controlled in Section~\ref{sec:analysis}.

\begin{remark}[Sharpness of the transport bound]
\label{rem:transport-sharpness}
Consider $p=3$, $n=1$, and $f(t):=t^4/24$. Then $\nabla^{p-1}f(t)=t^2/2$, $\nabla^{p}f(t)=t$, and Assumption~\ref{ass:main} holds with $L_3=1$. For arbitrary $D,E\geq0$, take $x:=0$, $y:=-E$, and $z:=D$, and suppose that $\widehat A_y=\nabla^{p-1}f(y)$ and $\widehat B_z=\nabla^{p}f(z)$. The affine transport rule gives
$\widehat A_x-\nabla^{p-1}f(x)=DE+\frac{E^2}{2}$.
Since $\norm{x-z}=D$ and $\norm{x-y}=E$, equality holds in \eqref{eq:transport}. Hence the coefficients of $L_p\norm{x-z}\norm{x-y}$ and $L_p\norm{x-y}^2/2$ are sharp for this affine transport rule.
If $D=0$, then $x=z$ and $\left|\widehat A_x-\nabla^{p-1}f(x)\right|=E^2/2$. Thus the quadratic term remains necessary when the current point coincides with the block reference. 
\end{remark}

\section{The multi-rate adaptive regularization method}
\label{sec:algorithm}

We now combine the finite-difference reconstructions and affine transport developed in Section~\ref{sec:reconstruction} to define the algorithm. The method proceeds in attempts. Each attempt uses a fixed estimate of the Lipschitz constant and generates a candidate block of at most $m$ model steps. The order-$p$ approximation is reconstructed once at the beginning of the block, whereas the order-$(p-1)$ approximation is reconstructed every $\ell$ steps and transported between successive reconstruction points, where $1\leq \ell\leq m$. Unless an $\eps$-approximate first-order stationary point is encountered earlier, the attempt is declared successful or unsuccessful after all $m$ steps have been generated.

\subsection{Parameters, reconstructions, and local model}

Let $L_k\geq L_{\min}>0$ denote the estimate of $L_p$ used during attempt $k$. Let $C_\sigma>0$ be a constant depending only on $p$,
specified in Appendix~\ref{app:constants}, and set
\begin{equation}\label{eq:sigma-k}
\sigma_k:=C_\sigma L_km\ell.
\end{equation}
The factor $m\ell$ reflects the two displacement scales in the transport bound \eqref{eq:transport} and controls the accumulated transport errors; see Section~\ref{sec:analysis}.

Set
\begin{equation}\label{eq:delta-bars}
 \ov\delta_{A,k}:=c_A\sigma_k^{(p-2)/p}\eps^{2/p},\qquad  \ov\delta_{B,k}:=c_B\ell^{-1}\sigma_k^{(p-1)/p}\eps^{1/p},
 \end{equation}
and choose
\begin{equation}\label{eq:hA-k}
h_{A,k}:=\left(\frac{6 \ov\delta_{A,k}}{\sqrt n\,L_k}\right)^{1/2},
\end{equation}
and
\begin{equation}\label{eq:hB-k}
h_{B,k} :=\frac{3 \ov\delta_{B,k}}{2nL_k}.
\end{equation}
The constants $c_A,c_B>0$, together with $\theta$ and $\eta_{\rm blk}$ used below, depend only on $p$ and are fixed as in Appendix~\ref{app:constants}. In particular, $0<\theta\leq\frac{1}{2}$.

Suppose that $L_k\geq L_p$ and let $x\in\R^n$ be any point at which a fresh reconstruction is formed during attempt $k$. Lemmas~\ref{lem:rec:p-1} and~\ref{lem:nab-p-reconstruction}, together with \eqref{eq:hA-k}--\eqref{eq:hB-k}, give
\begin{equation}\label{eq:fresh-A-error}
\left\|\sym\left(\widetilde A_{h_{A,k}}(x)\right)-\nabla^{p-1}f(x)\right\|_{[p-1]}\leq\frac{\sqrt n\,L_p}{6}h_{A,k}^2=\frac{L_p}{L_k} \ov\delta_{A,k}\leq \ov\delta_{A,k},
\end{equation}
and
\begin{equation}\label{eq:fresh-B-error}
 \left\|\sym\left(\widetilde B_{h_{B,k}}(x)\right)-\nabla^pf(x)\right\|_{[p]}\leq\frac{2nL_p}{3}h_{B,k}=\frac{L_p}{L_k} \ov\delta_{B,k}\leq \ov\delta_{B,k}.
\end{equation}

Let $z_k$ denote the accepted reference point at the beginning of attempt $k$, and set $x_{k,0}:=z_k$. The order-$p$ approximation used throughout the block is reconstructed once at $z_k$:
\begin{equation}\label{eq:B-k}
\widehat B_k:=\sym\left(\widetilde B_{h_{B,k}}(z_k)\right).
\end{equation}
For $t\in\{0,\ldots,m-1\}$, define the first index of the subblock containing $t$ by
\begin{equation}\label{eq:subblock-start}
b(t):=\ell\left\lfloor\frac{t}{\ell}\right\rfloor.
\end{equation}
Thus, $b(t)$ is the most recent subblock start at or before $t$, and the final subblock is automatically truncated if $\ell$ does not divide $m$. For each subblock start $b\in\{0,\ell,2\ell,\ldots,\ell\lfloor(m-1)/\ell\rfloor\}$, reconstruct
\begin{equation}\label{eq:A-kb}
\widehat A_{k,b}:=\sym\left(\widetilde A_{h_{A,k}}(x_{k,b})\right).
\end{equation}
At a local iteration $t$ in this subblock, transport this approximation from $x_{k,b(t)}$ to $x_{k,t}$ according to
\begin{equation}\label{eq:A-kt}
\widehat A_{k,t}:=\widehat A_{k,b(t)}+\widehat B_k[x_{k,t}-x_{k,b(t)}].
\end{equation}
Before the model at $x_{k,t}$ is formed, the bundle $\cF(x_{k,t})$ is available. Thus the derivatives through order $p-2$ are known exactly at the current point, while the two unavailable Taylor tensors are represented by $\widehat A_{k,t}$ and $\widehat B_k$. For a displacement $s\in\R^n$, define the regularized model
\begin{equation}\label{eq:model}
M_{k,t}(s):=f(x_{k,t})+\sum_{j=1}^{p-2}\frac{1}{j!}\nabla^j f(x_{k,t})[s]^j+\frac{1}{(p-1)!}\widehat A_{k,t}[s]^{p-1}+\frac{1}{p!}\widehat B_k[s]^p+\frac{\sigma_k}{(p+1)!}\|s\|^{p+1}.
\end{equation}
At local iteration $t$, the model is used to select a displacement $s_{k,t}$ from the current point $x_{k,t}$, after which the next trial point is $x_{k,t+1}:=x_{k,t}+s_{k,t}$. A step $s_{k,t}$ is called admissible if it produces strict model decrease and satisfies the relative first-order residual condition, that is,
\begin{equation}\label{eq:model-decrease}
M_{k,t}(s_{k,t})<M_{k,t}(0)=f(x_{k,t}),
\end{equation}
and
\begin{equation}\label{eq:model-residual}
\|\nabla_sM_{k,t}(s_{k,t})\|\leq\theta\frac{\sigma_k}{p!}\|s_{k,t}\|^p.
\end{equation}

\begin{lemma}[Existence of an admissible model step]\label{lem:admissible-step}
Suppose that $\nabla f(x_{k,t})\neq0$. Then $M_{k,t}$ admits an admissible step.
\end{lemma}
\begin{proof}
The terms of $M_{k,t}$ other than the regularizer have degree at most $p$. Since $\sigma_k>0$, we have $M_{k,t}(s)\to+\infty$ as $\|s\|\to\infty$. Hence $M_{k,t}$ is coercive and attains a global minimum. Since,
\begin{equation*}
\nabla_sM_{k,t}(0)=\nabla f(x_{k,t})\neq0,
\end{equation*}
the origin is not a minimizer. Therefore every global minimizer $s^\star$ is nonzero and satisfies
\begin{equation*}
M_{k,t}(s^\star)<M_{k,t}(0),\qquad \nabla_sM_{k,t}(s^\star)=0.
\end{equation*}
Hence $s^\star$ satisfies \eqref{eq:model-decrease}--\eqref{eq:model-residual}.
\end{proof}
The global minimizer in Lemma~\ref{lem:admissible-step} is used only to establish existence of an admissible step. The algorithm may use any step satisfying \eqref{eq:model-decrease}--\eqref{eq:model-residual}. Computing such a step requires no additional bundled-oracle query. A finite exact-arithmetic realization is given in Section~\ref{sec:complexity}.

\subsection{Block acceptance and oracle accounting}
If none of the stationarity tests succeeds during attempt $k$, the attempt is called \textit{complete and nonterminal}. It then generates $m$ model steps. 
Let $x_k^{\rm best}$ be the first point attaining
\begin{equation*}
f(x_k^{\rm best})=\mint{0\leq j\leq m}f(x_{k,j}),
\end{equation*}
Define
\begin{equation}\label{eq:required-decrease}
\Delta_k^{\rm req}:=\eta_{\rm blk}L_k^{-1/p}m^{(p-1)/p}\ell^{-1/p}\eps^{(p+1)/p}.
\end{equation}
A complete nonterminal attempt $k$ is successful if
\begin{equation}\label{eq:block-success}
f(z_k)-f(x_k^{\rm best})\geq\Delta_k^{\rm req}.
\end{equation}
Section~\ref{sec:analysis} shows that \eqref{eq:block-success} holds whenever $L_k\geq L_p$. On a successful attempt, $x_k^{\rm best}$ becomes the next reference point and the estimate $L_k$ is decreased. On an unsuccessful attempt, the reference point is unchanged and $L_k$ is doubled.
The complete method is given in Algorithm~\ref{alg:MRAR}.

\begin{algorithm}[!htbp]
\caption{$\MRAR(p,p-2)$}
\label{alg:MRAR}
\begin{algorithmic}[1]
\Require $x_0\in\R^n$, $\eps>0$, $p\geq3$, integers $1\leq\ell\leq m$, estimates $L_0\geq L_{\min}>0$, constants chosen as in Appendix~\ref{app:constants}, and a model-step routine returning a step satisfying \eqref{eq:model-decrease}--\eqref{eq:model-residual}.

\State $z_0\gets x_0$.
\State Evaluate and retain $\cF(z_0)$.
\If{$\|\nabla f(z_0)\|\leq\eps$}
    \State \Return $z_0$.
\EndIf

\For{$k=0,1,2,\ldots$}
    \State $x_{k,0}\gets z_k$ and $x_k^{\rm best}\gets z_k$.
    \State Compute $\sigma_k$, $ \ov\delta_{A,k}$, $ \ov\delta_{B,k}$, $h_{A,k}$, and $h_{B,k}$ from \eqref{eq:sigma-k}--\eqref{eq:hB-k}.
    \State Reconstruct $\widehat B_k \gets \sym\left(\widetilde B_{h_{B,k}}(z_k)\right)$.

    \For{$t=0,\ldots,m-1$}

        \State $b\gets\ell\left\lfloor t/\ell\right\rfloor$.

        \If{$t=b$}
            \State Reconstruct $\widehat A_{k,b}\gets \sym\left(\widetilde A_{h_{A,k}}(x_{k,b})\right)$.
        \EndIf

        \State $\widehat A_{k,t}\gets\widehat A_{k,b}+\widehat B_k[x_{k,t}-x_{k,b}]$.
        \State Form $M_{k,t}$ from \eqref{eq:model}.
        \State Use the model-step routine to obtain an admissible step $s_{k,t}$.
        \State $x_{k,t+1}\gets x_{k,t}+s_{k,t}$.
        \State Evaluate $\cF(x_{k,t+1})$.

        \If{$f(x_{k,t+1})<f(x_k^{\rm best})$}
            \State $x_k^{\rm best}\gets x_{k,t+1}$ and retain its bundle.
        \EndIf

        \If{$\|\nabla f(x_{k,t+1})\|\leq\eps$}
            \State \Return $x_{k,t+1}$.
        \EndIf

    \EndFor

    \State Compute $\Delta_k^{\rm req}$ from \eqref{eq:required-decrease}.

    \If{$f(z_k)-f(x_k^{\rm best})\geq\Delta_k^{\rm req}$}
        \State $z_{k+1}\gets x_k^{\rm best}$ and retain $\cF(z_{k+1})$.
        \State $L_{k+1}\gets \max\{L_{\min},L_k/2\}$.
    \Else
        \State $z_{k+1}\gets z_k$ and retain $\cF(z_k)$.
        \State $L_{k+1}\gets 2L_k$.
    \EndIf

\EndFor

\end{algorithmic}
\end{algorithm}

Every point at which a model is formed has already failed the stationarity test, so $\|\nabla f(x_{k,t})\|>\eps$ and Lemma~\ref{lem:admissible-step} ensures that an admissible step exists.
At the beginning of each attempt, $\cF(z_k)$ is retained, and before $M_{k,t}$ is formed, $\cF(x_{k,t})$ is available. The bundle associated with $x_k^{\rm best}$ is also retained. Thus affine transport, best-point selection, and the model-step routine require no additional bundled-oracle queries. The exact query count is given in
Section~\ref{sec:complexity}.

The best-point rule allows individual local steps to increase the objective value. Only the decrease over the complete block is used in the acceptance test \eqref{eq:block-success}.

\section{Model errors and convergence analysis}
\label{sec:analysis}
This section establishes the convergence properties of Algorithm~\ref{alg:MRAR}. We first show that the reconstruction and reuse errors remain controlled after summation over a complete block. We then prove that every complete nonterminal attempt is successful whenever $L_k\geq L_p$, which yields a uniform bound on the adaptive estimates and finite termination.

\subsection{Local model estimates and blockwise error control}

Fix a complete nonterminal attempt $k$ with $L_k\geq L_p$. 
By the choice of $C_\sigma$ in Appendix~\ref{app:constants}, $C_\sigma\geq2$ and
\begin{equation}\label{eq:sigma-large}
\sigma_k=C_\sigma L_km\ell\geq 2L_p.
\end{equation}
For $t=0,\ldots,m-1$, define
\begin{equation}\label{eq:local-errors}
\delta_{A,k,t}:=\|\widehat A_{k,t}-\nabla^{p-1}f(x_{k,t})\|_{[p-1]},\qquad \delta_{B,k,t}:=\|\widehat B_k-\nabla^pf(x_{k,t})\|_{[p]}.
\end{equation}
These quantities include both the reconstruction error and the error caused by transport or reuse away from the reconstruction point. Set
\begin{equation}\label{eq:local-error-penalty}
\mathcal E_{k,t}:=\frac{\delta_{A,k,t}^{(p+1)/2}}{\sigma_k^{(p-1)/2}}+\frac{\delta_{B,k,t}^{p+1}}{\sigma_k^p}.
\end{equation}
The polynomial part of $M_{k,t}$ in \eqref{eq:model} differs from the exact degree-$p$ Taylor polynomial $T_p(x_{k,t},s)$ only in derivative orders $p-1$ and $p$. Hence, by \eqref{eq:function-rem}--\eqref{eq:gradient-rem} and \eqref{eq:local-errors}, for every $s\in\R^n$,
\begin{equation}\label{eq:approx-value-error}
\left|f(x_{k,t}+s)-M_{k,t}(s)+\frac{\sigma_k}{(p+1)!}\|s\|^{p+1}\right|\leq\frac{L_p}{(p+1)!}\|s\|^{p+1}
+\frac{\delta_{A,k,t}}{(p-1)!}\|s\|^{p-1}+\frac{\delta_{B,k,t}}{p!}\|s\|^p,
\end{equation}
and
\begin{equation}\label{eq:approx-gradient-error}
\left\|\nabla f(x_{k,t}+s)-\nabla_sM_{k,t}(s)+\frac{\sigma_k}{p!}\|s\|^{p-1}s\right\|
\leq\frac{L_p}{p!}\|s\|^p+\frac{\delta_{A,k,t}}{(p-2)!}\|s\|^{p-2}+\frac{\delta_{B,k,t}}{(p-1)!}\|s\|^{p-1}.
\end{equation}
The admissibility conditions and \eqref{eq:approx-value-error}--\eqref{eq:approx-gradient-error} give the following local estimates.

\begin{lemma}[Local step estimates]\label{lem:local-step-estimates}
Suppose that $L_k\geq L_p$.
Then, for every $t=0,\ldots,m-1$,
\begin{equation}\label{eq:local-decrease}
f(x_{k,t})-f(x_{k,t+1})\geq c_f\sigma_k\|s_{k,t}\|^{p+1}-K_{f,p}\mathcal E_{k,t},
\end{equation}
and
\begin{equation}\label{eq:local-gradient}
\|\nabla f(x_{k,t+1})\|^{(p+1)/p}\leq C_{g,p}\sigma_k^{1/p}\left(\sigma_k\|s_{k,t}\|^{p+1}+\mathcal E_{k,t}\right).
\end{equation}
The constants $c_f$, $K_{f,p}$ and $C_{g,p}$ depend only on $p$ and are specified in Appendix~\ref{app:constants}.
\end{lemma}
\begin{proof}
See Appendix~\ref{sec:app:proofs}.
\end{proof}

Since the attempt is complete and nonterminal, $\|\nabla f(x_{k,t+1})\|>\eps$ for every $t=0,\ldots,m-1$. Hence \eqref{eq:local-gradient} gives
\begin{equation}\label{eq:nonterminal-local-gradient}
\eps^{(p+1)/p}<C_{g,p}\sigma_k^{1/p}\left(\sigma_k\|s_{k,t}\|^{p+1}+\mathcal E_{k,t}\right),\qquad t=0,\ldots,m-1.
\end{equation}
We next control $\mathcal E_{k,t}$ over the complete block.  Adding and subtracting $\nabla^pf(z_k)$ gives
\[
\widehat B_k-\nabla^pf(x_{k,t})=\left(\widehat B_k-\nabla^pf(z_k)\right)+\left(\nabla^pf(z_k)-\nabla^pf(x_{k,t})\right),
\]
Since $\widehat B_k$ is reconstructed at $z_k$, \eqref{eq:fresh-B-error} applies to the first term. Hence, using also \eqref{eq:Lp}, 
\begin{equation}\label{eq:B-error-along-block}
\delta_{B,k,t}\leq\ov\delta_{B,k}+L_p\|x_{k,t}-z_k\|.
\end{equation}
For the order-$(p-1)$ approximation, apply \eqref{eq:transport} with $z=z_k$, $y=x_{k,b(t)}$, and $x=x_{k,t}$.
Using
\[
\|\widehat A_{k,b(t)}-\nabla^{p-1}f(x_{k,b(t)})\|_{[p-1]}\leq\ov\delta_{A,k},\qquad
\|\widehat B_k-\nabla^pf(z_k)\|_{[p]}\leq\ov\delta_{B,k},
\]
we obtain
\begin{equation}\label{eq:A-error-along-block}
\delta_{A,k,t}\leq\ov\delta_{A,k}+\ov\delta_{B,k}\|x_{k,t}-x_{k,b(t)}\|+L_p\|x_{k,t}-z_k\|\|x_{k,t}-x_{k,b(t)}\|+\frac{L_p}{2}\|x_{k,t}-x_{k,b(t)}\|^2.
\end{equation}
For the complete block, set
\begin{equation}\label{eq:block-energy}
S_k:=\sum_{t=0}^{m-1}\Vert s_{k,t}\Vert^{p+1}.
\end{equation}
The displacement terms in \eqref{eq:B-error-along-block}--\eqref{eq:A-error-along-block} are controlled by the two-scale path estimates in Lemma~\ref{lem:path-bounds}.

\begin{lemma}[Accumulated reconstruction and reuse error]\label{lem:accumulated-error}
For every complete nonterminal attempt $k$ with $L_k\geq L_p$,
\begin{equation}\label{eq:accumulated-error}
\sum_{t=0}^{m-1}\mathcal E_{k,t}\leq\rho\sigma_k S_k+\eta_{\rm err}m\eps^{(p+1)/p}\sigma_k^{-1/p}.
\end{equation}
\end{lemma}
\begin{proof}
See Appendix~\ref{sec:app:proofs}.
\end{proof}
The first term in \eqref{eq:accumulated-error} controls the errors caused by transport and reuse, while the second contains the remaining fresh-reconstruction errors.

\subsection{Successful blocks and finite termination}

We first show that every complete nonterminal attempt is successful whenever $L_k\geq L_p$.

\begin{proposition}[Success for sufficiently large estimates]
\label{prop:block-success}
Let attempt $k$ be complete and nonterminal. If $L_k\geq L_p$, then
\begin{equation}\label{eq:large-estimate-success}
f(z_k)-f(x_k^{\rm best})>\Delta_k^{\rm req}.
\end{equation}
In particular, attempt $k$ is successful.
\end{proposition}
\begin{proof}
Summing \eqref{eq:nonterminal-local-gradient} over $t=0,\ldots,m-1$ and applying Lemma~\ref{lem:accumulated-error} gives
\begin{equation*}
m\eps^{(p+1)/p}<C_{g,p}\sigma_k^{1/p}\left[(1+\rho)\sigma_kS_k+\eta_{\rm err}m\eps^{(p+1)/p}\sigma_k^{-1/p}\right].
\end{equation*}
Since $C_{g,p}\eta_{\rm err}\leq1/4$, it follows that
\begin{equation}\label{eq:block-step-lower}
\sigma_kS_k\geq\frac{3}{4C_{g,p}(1+\rho)}m\eps^{(p+1)/p}\sigma_k^{-1/p}.
\end{equation}
Summing \eqref{eq:local-decrease} over the complete block and applying Lemma~\ref{lem:accumulated-error} again yields
\begin{equation}\label{eq:block-decrease-before-constants}
f(z_k)-f(x_{k,m})\geq(c_f-K_{f,p}\rho)\sigma_kS_k-K_{f,p}\eta_{\rm err}m\eps^{(p+1)/p}\sigma_k^{-1/p}.
\end{equation}
The choices of $\rho$ and $\eta_{\rm err}$ in Appendix~\ref{app:constants}, together with \eqref{eq:block-step-lower}, imply
\[
f(z_k)-f(x_{k,m}) \geq \frac{3c_fc_E}{8} m\eps^{(p+1)/p}\sigma_k^{-1/p} = \eta_0 m\eps^{(p+1)/p}\sigma_k^{-1/p},
\]
where $\eta_0$ is specified in Appendix~\ref{app:constants}. Since $f(x_k^{\rm best})\leq f(x_{k,m})$ and $\sigma_k=C_\sigma L_km\ell$, we obtain
\begin{align*}
f(z_k)-f(x_k^{\rm best}) &\geq \eta_0C_\sigma^{-1/p} L_k^{-1/p} m^{(p-1)/p} \ell^{-1/p} \eps^{(p+1)/p}\\
&= \eta_{\rm th} L_k^{-1/p} m^{(p-1)/p} \ell^{-1/p} \eps^{(p+1)/p},
\end{align*}
where $\eta_{\rm th}:=\eta_0C_\sigma^{-1/p}$. Since $\eta_{\rm blk}<\eta_{\rm th}$, the definition \eqref{eq:required-decrease} gives
\[
f(z_k)-f(x_k^{\rm best})> \eta_{\rm blk} L_k^{-1/p} m^{(p-1)/p} \ell^{-1/p} \eps^{(p+1)/p} = \Delta_k^{\rm req}.
\]
This proves \eqref{eq:large-estimate-success}.
\end{proof}

By Proposition~\ref{prop:block-success}, a complete nonterminal attempt can be unsuccessful only if $L_k<L_p$.
\begin{theorem}[Finite termination and attempt bound]\label{thm:termination}
Let $L_{\max}:=\max\{L_0,2L_p\}$ and $\Delta_f:=f(x_0)-f_{\rm low}$. Then, 
\begin{enumerate}[label=(\textbf{\alph*}), font=\normalfont\bfseries, leftmargin=0.7cm]
\item\label{thm:termination:a} for every attempt $k$,
\begin{equation}\label{eq:Lk-bound}
L_k\leq L_{\max}.
\end{equation}
 
\item\label{thm:termination:b} if $N_{\rm succ}$ and $N_{\rm fail}$ denote, respectively, the numbers of completed successful and unsuccessful nonterminal attempts before termination, then
\begin{equation}\label{eq:successful-attempts}
N_{\rm succ}\leq\frac{L_{\max}^{1/p}\Delta_f}{\eta_{\rm blk}}m^{-(p-1)/p}\ell^{1/p}\eps^{-(p+1)/p},
\end{equation}
and
\begin{equation}\label{eq:failed-attempts}
N_{\rm fail}\leq N_{\rm succ}+\log_2\left(\frac{L_{\max}}{L_0}\right).
\end{equation}

\item\label{thm:termination:c} with $N_{\rm comp}:=N_{\rm succ}+N_{\rm fail}$, we have
\begin{equation}\label{eq:completed-attempts}
N_{\rm comp}\leq\frac{2L_{\max}^{1/p}\Delta_f}{\eta_{\rm blk}}m^{-(p-1)/p}\ell^{1/p}\eps^{-(p+1)/p}+\log_2\left(\frac{L_{\max}}{L_0}\right).
\end{equation}

\item\label{thm:termination:d} Algorithm~\ref{alg:MRAR} terminates after finitely many attempts at a point $x$ satisfying $\|\nabla f(x)\|\leq\eps$. 
\end{enumerate}
\end{theorem}
\begin{proof}
\ref{thm:termination:a} The claim is immediate for $k=0$. Suppose that attempt $k$ begins with $L_k\leq L_{\max}$ and does not terminate early. If the attempt is successful, then
\[
L_{k+1}=\max\{L_{\min},L_k/2\}\leq L_k\leq L_{\max}.
\]
If the attempt is unsuccessful, Proposition~\ref{prop:block-success} implies $L_k<L_p$, and therefore
\[
L_{k+1}=2L_k<2L_p\leq L_{\max}.
\]
Hence \eqref{eq:Lk-bound} follows by induction.
\\
\ref{thm:termination:b} By the acceptance condition and \eqref{eq:Lk-bound}, every completed successful attempt satisfies
\[
f(z_k)-f(z_{k+1})\geq\eta_{\rm blk}L_k^{-1/p}m^{(p-1)/p}\ell^{-1/p}\eps^{(p+1)/p}\geq\eta_{\rm blk}L_{\max}^{-1/p}m^{(p-1)/p}\ell^{-1/p}\eps^{(p+1)/p}.
\]
Failed attempts leave the accepted reference unchanged. Since $f\geq f_{\rm low}$, summing the decreases over all completed successful attempts gives
\[
N_{\rm succ}\eta_{\rm blk}L_{\max}^{-1/p}m^{(p-1)/p}\ell^{-1/p}\eps^{(p+1)/p}\leq f(x_0)-f_{\rm low}=\Delta_f.
\]
Rearranging proves \eqref{eq:successful-attempts}.
It remains to count unsuccessful attempts. Define $r_k:=\log_2\left(\frac{L_k}{L_0}\right)$.
Every unsuccessful attempt doubles $L_k$ and therefore increases $r_k$ by exactly one. On a successful attempt,
\[
L_{k+1}=\max\{L_{\min},L_k/2\}\geq L_k/2,
\]
so $r_{k+1}-r_k\geq-1$. Consequently, after any finite sequence of completed nonterminal attempts,
\[
N_{\rm fail}-N_{\rm succ}\leq r_k=\log_2\left(\frac{L_k}{L_0}\right).
\]
Using \eqref{eq:Lk-bound},
\[
N_{\rm fail}-N_{\rm succ}\leq\log_2\left(\frac{L_{\max}}{L_0}\right),
\]
which proves \eqref{eq:failed-attempts}.
\\
\ref{thm:termination:c} Combining \eqref{eq:successful-attempts} and \eqref{eq:failed-attempts} yields \eqref{eq:completed-attempts}.
\\
\ref{thm:termination:d} Suppose that the algorithm does not terminate. Then every attempt is nonterminal and, since each attempt contains at most $m$ local iterations, every attempt is completed. An infinite run would therefore contain arbitrarily many completed nonterminal attempts, contradicting the finite bound \eqref{eq:completed-attempts}. Hence a stationarity test succeeds after finitely many attempts, and Algorithm~\ref{alg:MRAR} terminates at a point $x$ satisfying $\|\nabla f(x)\|\leq\eps$.
\end{proof}

\section{Oracle complexity and computational cost}
\label{sec:complexity}

We now translate the attempt bound of Section~\ref{sec:analysis} into a count of bundled order-$(p-2)$ oracle queries. We first give the exact point-query cost for arbitrary block and subblock lengths and derive a general complexity bound. We then choose the two refresh lengths to obtain the dimension dependence of the accuracy-dependent term.

\subsection{Exact point-query count and general complexity}
For a complete nonterminal attempt, the number of new bundled order-$(p-2)$ point queries is
\begin{equation}\label{eq:Ccomplete}
C_{\rm comp}:=n(n+1)+2n\left\lceil\frac{m}{\ell}\right\rceil+m.
\end{equation}
The three terms correspond, respectively, to the order-$p$ reconstruction, the order-$(p-1)$ reconstructions, and the $m$ trial points.

If an attempt terminates after generating $j\in\{1,\ldots,m\}$ trial points, then its number of new queries is
\begin{equation}\label{eq:Cterminal}
C_{\rm term}(j):=n(n+1)+2n\left\lceil\frac{j}{\ell}\right\rceil+j.
\end{equation}
Indeed, the executed local iterations are $t=0,\ldots,j-1$, so exactly $\lceil j/\ell\rceil$ subblocks have started. Hence
\begin{equation}\label{eq:Cterminal-bound}
C_{\rm term}(j)\leq C_{\rm comp}.
\end{equation}

\begin{lemma}[Point-query count]
\label{lem:oracle-count}
If $\|\nabla f(x_0)\|\leq\eps$, then $N_{p-2}(\eps)=1$. Otherwise, let $N_{\rm comp}$ be the number of completed nonterminal attempts before termination, and suppose that the terminal attempt generates $j_\star\in\{1,\ldots,m\}$ trial points. Then
\begin{equation}\label{eq:exact-ledger}
N_{p-2}(\eps)=1+C_{\rm comp}N_{\rm comp}+C_{\rm term}(j_\star).
\end{equation}
\end{lemma}
\begin{proof}
The initial bundle $\cF(x_0)$ contributes one query. Every completed nonterminal attempt contributes $C_{\rm comp}$ new queries, while the terminal attempt contributes $C_{\rm term}(j_\star)$ queries. The reference bundle at the beginning of each attempt is already available by construction of Algorithm~\ref{alg:MRAR}. This proves \eqref{eq:exact-ledger}.
\end{proof}

\begin{theorem}[Bundled order-$(p-2)$ oracle complexity]
\label{thm:main}
Let Assumption~\ref{ass:main} hold, choose the algorithmic constants as in Appendix~\ref{app:constants}, and fix integers
$1\leq\ell\leq m$. Then Algorithm~\ref{alg:MRAR} terminates at a point
$x$ satisfying $\|\nabla f(x)\|\leq\eps$, and
\begin{equation}\label{eq:main-bound}
\boxed{N_{p-2}(\eps)\leq1+C_{\rm comp}\left[\frac{2L_{\max}^{1/p}\Delta_f}{\eta_{\rm blk}}m^{-(p-1)/p}\ell^{1/p}\eps^{-(p+1)/p}+\log_2\left(\frac{L_{\max}}{L_0}\right)+1\right],}
\end{equation}
where
$\Delta_f:=f(x_0)-f_{\rm low}$ and $L_{\max}:=\max\{L_0,2L_p\}$.
\end{theorem}
\begin{proof}
Finite termination follows from Theorem~\ref{thm:termination}. If $\|\nabla f(x_0)\|\leq\eps$, then $N_{p-2}(\eps)=1$. Otherwise, Lemma~\ref{lem:oracle-count} and \eqref{eq:Cterminal-bound} give
\[
N_{p-2}(\eps)\leq 1+C_{\rm comp}(N_{\rm comp}+1).
\]
The result follows from \eqref{eq:completed-attempts}.
\end{proof}

\subsection{Choice of the block and subblock lengths}

The accuracy-dependent term in \eqref{eq:main-bound} is proportional to
\begin{equation}\label{eq:Phi-definition}
\Phi_n(m,\ell):=C_{\rm comp}m^{-(p-1)/p}\ell^{1/p} = \left[n(n+1)+2n\left\lceil\frac{m}{\ell}\right\rceil+m\right]m^{-(p-1)/p}\ell^{1/p}.
\end{equation}
We next determine the dimension dependence obtained by choosing $m$ and $\ell$.

\begin{proposition}[Two-rate block-length scaling]\label{prop:architecture-optimality}
For every $n\geq1$,
\begin{equation}\label{eq:architecture-optimality}
\inft{\substack{m,\ell\in\N\\1\leq\ell\leq m}} \Phi_n(m,\ell) = \Theta_p(n^{3/p}).
\end{equation}
\end{proposition}
\begin{proof}
Using $n(n+1)\geq n^2$ and $\lceil m/\ell\rceil\geq m/\ell$, we obtain
\[
\Phi_n(m,\ell) \geq \left(n^2+\frac{2nm}{\ell}+m\right)m^{-(p-1)/p}\ell^{1/p}.
\]
Set $m=n^2u$ and $\ell=nv$, with $u,v>0$. Then $\Phi_n(m,\ell)\geq n^{3/p}\psi_p(u,v)$, where
\begin{equation}\label{eq:psi-definition}
\psi_p(u,v) := u^{-(p-1)/p}v^{1/p} + u^{1/p}v^{1/p} + 2u^{1/p}v^{-(p-1)/p}.
\end{equation}
Let
\[
a:=u^{-(p-1)/p}v^{1/p},
\qquad
b:=u^{1/p}v^{1/p},
\qquad
c:=2u^{1/p}v^{-(p-1)/p}.
\]
Since $ab^{p-2}c=2$, the weighted arithmetic--geometric mean inequality gives
\begin{align*}
\psi_p(u,v) = \frac1p(pa) + \frac{p-2}{p} \left(\frac{p}{p-2}b\right) + \frac1p(pc)
\geq p\,2^{1/p}(p-2)^{-(p-2)/p}. 
\end{align*}
Thus, there exists $c_p>0$, depending only on $p$, such that
\[
\inft{\substack{m,\ell\in\N\\1\leq\ell\leq m}}\Phi_n(m,\ell)\geq c_p n^{3/p}.
\]
For the reverse bound, choose $m=(p-2)n^2$ and $\ell=(p-2)n$. Then $m/\ell=n$ and
\[
C_{\rm comp} = n(n+1)+2n^2+(p-2)n^2 = (p+1)n^2+n. 
\]
Therefore,
\[
\Phi_n(m,\ell) = (p-2)^{-(p-2)/p}\left(p+1+\frac1n \right)n^{3/p}.
\]
Hence there exists $C_p>0$, depending only on $p$, such that
\[
\inft{\substack{m,\ell\in\N\\1\leq\ell\leq m}} \Phi_n(m,\ell)\leq C_p n^{3/p}.
\]
Combining the two bounds proves \eqref{eq:architecture-optimality}.
\end{proof}

Motivated by Proposition~\ref{prop:architecture-optimality}, we use the simple integer choice
\begin{equation}\label{eq:canonical}
m:=(p-2)n^2,\qquad \ell:=(p-2)n.
\end{equation}
This choice is feasible for every $p\geq3$ and $n\geq1$, and satisfies $m/\ell=n$. Hence
\begin{equation}\label{eq:canonical-Ccomp}
C_{\rm comp}=(p+1)n^2+n,
\end{equation}
and
\begin{equation}\label{eq:canonical-leading-factor}
C_{\rm comp}m^{-(p-1)/p}\ell^{1/p}=(p-2)^{-(p-2)/p}\left(p+1+\frac1n\right)n^{3/p}.
\end{equation}

\begin{corollary}[Canonical point-query bound]\label{cor:canonical}
Under the assumptions of Theorem~\ref{thm:main}, the choice \eqref{eq:canonical} gives
\begin{equation}\label{eq:canonical-bound}
\boxed{N_{p-2}(\eps)=\mathcal O_p\left(n^{3/p}L_{\max}^{1/p}\Delta_f\eps^{-(p+1)/p}+n^2\left[1+\log_2\left(\frac{L_{\max}}{L_0}\right)\right]\right).}
\end{equation}
\end{corollary}

\begin{remark}[Benefit of separate refresh lengths]
\label{rem:architecture-scope}
The positive relaxation used in the lower-bound proof of Proposition~\ref{prop:architecture-optimality} attains equality at
$u=p-2$ and $v=2(p-2)$. The integer choice \eqref{eq:canonical} is feasible for every $n\geq1$ and attains the same $\Theta_p(n^{3/p})$ scaling. The benefit of using separate block and subblock lengths can also be
seen directly from \eqref{eq:Phi-definition}. Reconstructing both
missing tensors at every step gives
\[
\Phi_n(1,1)=\Theta(n^2).
\]
If a single refresh length is imposed by setting $\ell=m$, then
\[
\Phi_n(m,m) = \left[n(n+1)+2n+m\right]m^{-(p-2)/p}.
\]
The continuous relaxation is minimized at
\[
m_\star=\frac{p-2}{2}(n^2+3n),
\]
and hence
\[
\inft{m>0}\Phi_n(m,m)= \Theta_p(n^{4/p}).
\]
For the integer problem, choose $\ov m:=\lceil m_\star\rceil$. Since $m_\star=\Theta_p(n^2)$,
$\Phi_n(\ov m,\ov m)=\mathcal O_p(n^{4/p})$.
Together with the lower bound provided by the continuous relaxation, this gives
\[
\inft{m\in\N}\Phi_n(m,m)=\Theta_p(n^{4/p}).
\]
Thus allowing separate block and subblock lengths improves the dimension factor from $\Theta_p(n^{4/p})$ to
$\Theta_p(n^{3/p})$.
\end{remark}

For $p=3$, \eqref{eq:canonical-bound} gives the accuracy-dependent term $n\eps^{-4/3}$. This is consistent with the dimension-free first-order lower bound of order $\eps^{-5/3}$ under Lipschitz third derivatives in \cite{Zhou2026FirstOrder}, since that lower bound is taken over a function class containing problems of arbitrary dimension, whereas \eqref{eq:canonical-bound} tracks $n$ and $\eps$ jointly.

\subsection{Inner model step and computational cost}

The point-query complexity above does not include the work required to compute an admissible model step. We first show that the conditions \eqref{eq:model-decrease}--\eqref{eq:model-residual} can be satisfied after finitely many iterations using only the stored model.

Fix a model $M$ of the form \eqref{eq:model} at a nonterminal point, so that $\nabla_sM(0)\neq0$, and let $\sigma>0$ denote its regularization parameter. Set $s_0:=0$ and choose $\alpha_{\max}>0$ and $\beta,c\in(0,1)$. At iteration $j$, let
\[
\alpha_j:=\alpha_{\max}\beta^r,
\]
where $r\in\Nz$ is the smallest integer such that
\begin{equation}\label{eq:inner-armijo}
M(s_j-\alpha_j\nabla_sM(s_j))\leq M(s_j)-c\alpha_j\|\nabla_sM(s_j)\|^2.
\end{equation}
and set
\begin{equation}\label{eq:inner-update}
s_{j+1}:=s_j-\alpha_j\nabla_sM(s_j).
\end{equation}

\begin{proposition}[Finite realization of an admissible step]\label{prop:finite-inner-realization}
The Armijo iteration \eqref{eq:inner-armijo}--\eqref{eq:inner-update} reaches, after finitely many accepted steps, an iterate $s_j$ satisfying
\begin{equation}\label{eq:inner-admissible}
M(s_j)<M(0),\qquad \|\nabla_sM(s_j)\|\leq\theta\frac{\sigma}{p!}\|s_j\|^p.
\end{equation}
Hence $s_j$ is an admissible model step.
\end{proposition}
\begin{proof}
See Appendix~\ref{sec:app:proofs}.
\end{proof}

Proposition~\ref{prop:finite-inner-realization} shows that the admissible-step condition can be realized using only evaluations of the stored model and its gradient, and therefore requires no additional queries to the original oracle.

The point-query measure also does not account for the amount of information returned by each oracle query. A bundled order-$(p-2)$ response contains
\[
1+\sum_{j=1}^{p-2}\binom{n+j-1}{j} = \binom{n+p-2}{p-2}= \Theta_p(n^{p-2}),
\]
independent scalar components. Similarly, a dense symmetric order-$p$ tensor has
\[
\binom{n+p-1}{p} = \Theta_p(n^p),
\]
independent entries. Thus the point-query bound in Theorem~\ref{thm:main} should not be interpreted as an arithmetic or memory complexity bound.

The cost of solving the regularized polynomial subproblem can also be substantial. For $p=3$, specialized methods for quartically regularized cubic models are developed in \cite{CartisZhu2025,CartisZhu2026}, while practical implementations and subproblem-selection strategies for third-order adaptive regularization are studied in \cite{CartisEtAl2026AR3}. Finite-precision, matrix-free, and uniform inner-complexity guarantees are beyond the scope of this paper.

\section{Conclusions}
\label{sec:conclusion}

We considered high-order adaptive regularization for functions with Lipschitz continuous $p$th derivatives, $p\geq3$, when derivative information is available only through order $p-2$. The two unavailable Taylor tensors are reconstructed from finite differences of $\nabla^{p-2}f$. The order-$p$ approximation is reconstructed once per block, while the order-$(p-1)$ approximation is reconstructed on subblocks and transported between reconstruction points.

For $m=(p-2)n^2$ and $\ell=(p-2)n$, the method finds an $\eps$-approximate first-order stationary point using
\[
\mathcal O_p\left(n^{3/p}L_{\max}^{1/p}\Delta_f\eps^{-(p+1)/p}+n^2\left[1+\log_2\left(\frac{L_{\max}}{L_0}\right)\right]
\right).
\]
bundled order-$(p-2)$ point queries. Thus the accuracy-dependent term retains the $\eps^{-(p+1)/p}$ dependence of full $p$th-order regularization. For the dense reconstruction considered here, the use of separate block and subblock lengths gives the dimension factor $\Theta_p(n^{3/p})$, compared with $\Theta_p(n^{4/p})$ when a common length is imposed.

The present analysis concerns dense finite-difference reconstructions and point-query complexity. Natural directions for further work include structured or matrix-free reconstructions and lower-order implementations based on oracles of order below $p-2$.


\section*{Statements and Declarations}

\subsection*{Funding}
The author was supported by the Research Foundation--Flanders (FWO) through research project G081222N and by the University of Antwerp Research Fund (BOF) through DOCPRO4 project 46929.

\subsection*{Competing interests}
The author has no relevant financial or non-financial interests to disclose.

\subsection*{Data availability}
No datasets were generated or analyzed in this theoretical study.

\appendix

\section{Choice of constants}\label{app:constants}

All constants below are chosen independently of $n$, $\eps$, and the problem data, and may depend on $p$. We specify them in the order in which they are needed.

Define
\[
R:=\frac{p+1}{2},
\qquad
\tau:=\frac{p+1}{p},
\qquad
H_p:=4^{(p-1)/2},
\qquad
c_f:=\frac{1}{4(p+1)!}.
\]
Set
\begin{align*}
K_{A,p}:=\frac{[8p(p-1)]^{(p-1)/2}}{R(p-1)!},
\qquad
K_{B,p}:=\frac{(8p)^p}{(p+1)p!},\qquad
K_{f,p}:=\max\{K_{A,p},K_{B,p}\}.
\end{align*}
Choose any fixed $0<\theta\leq1/2$ and define
\begin{align*}
    &a_0:=\frac{\theta+2}{p!},
    \qquad a_A:=\frac1{(p-2)!},
    \qquad a_B :=\frac1{(p-1)!},
    \\
  &  C_{g,p}:= 3^{1/p}\max\left\{a_0^\tau+\frac{p-2}{p}a_A^\tau+\frac{p-1}{p}a_B^\tau,
    \frac{2}{p}a_A^\tau, \frac{1}{p}a_B^\tau\right\}.
\end{align*}
Let
\[
\rho:=\min\left\{\frac14,\frac{c_f}{4K_{f,p}}\right\}.
\]
Then
\[
c_f-K_{f,p}\rho\geq\frac{3c_f}{4}.
\]
Choose $C_\sigma>0$ such that
\begin{equation}\label{eq:Csigma-choice}
C_\sigma
\geq
\max\left\{
2,
\left(\frac{4H_p}{(p+1)\rho}\right)^{1/R},
\left(\frac{2^{p+2}}{(p+1)\rho}\right)^{1/(p+1)}
\right\}.
\end{equation}
The first condition in \eqref{eq:Csigma-choice} ensures $\sigma_k\geq2L_p$ whenever $L_k\geq L_p$. The remaining conditions provide the relative reserves used to absorb the path-dependent reconstruction and transport errors in the proof of Lemma~\ref{lem:accumulated-error}.

Define
\[
c_E:=\frac{3}{4C_{g,p}(1+\rho)},
\]
and
\[
\eta_{\rm err}:=\min\left\{\frac{1}{4C_{g,p}},\frac{3c_fc_E}{8K_{f,p}}\right\}.
\]
In particular,
\[
C_{g,p}\eta_{\rm err}\leq\frac14,\qquad K_{f,p}\eta_{\rm err}\leq\frac{3c_fc_E}{8}.
\]
Choose $c_A,c_B>0$ such that
\begin{equation}\label{eq:fresh-error-choices}
H_pc_A^R\leq\frac{\eta_{\rm err}}{2}, 
\qquad
\left(2^p+\frac{H_p^2}{(p+1)\rho}\right)c_B^{p+1}\leq \frac{\eta_{\rm err}}{2}.
\end{equation}
These inequalities allocate the remaining fresh-reconstruction error in Lemma~\ref{lem:accumulated-error} between the order-$(p-1)$ and order-$p$ reconstructions.

Finally, set
\[
\eta_0:=\frac{3c_fc_E}{8},
\qquad
\eta_{\rm th}:=\eta_0C_\sigma^{-1/p},
\]
and choose
\[
0<\eta_{\rm blk}<\eta_{\rm th}.
\]

\section{Auxiliary results}
\label{app:auxiliary}
\renewcommand{\thetheorem}{B.\arabic{theorem}}
\setcounter{theorem}{0} 
The following estimates control the displacements generated within a complete block and within its subblocks.

\begin{lemma}[Two-scale path bounds]\label{lem:path-bounds}
Let $k$ be a complete nonterminal attempt, and let $S_k$ be defined by \eqref{eq:block-energy}. Then
\begin{align}
&\sum_{t=0}^{m-1}\|x_{k,t}-z_k\|^{p+1}\leq\frac{m^{p+1}}{p+1}S_k,\label{eq:block-path}\\
&\sum_{t=0}^{m-1}\|x_{k,t}-x_{k,b(t)}\|^{p+1}\leq\frac{\ell^{p+1}}{p+1}S_k,\label{eq:subblock-path}\\
&\sum_{t=0}^{m-1}\left(\|x_{k,t}-z_k\|\|x_{k,t}-x_{k,b(t)}\|\right)^{(p+1)/2}\leq\frac{(m\ell)^{(p+1)/2}}{p+1}S_k,\label{eq:mixed-path}\\
&\sum_{t=0}^{m-1}\|x_{k,t}-x_{k,b(t)}\|^{(p+1)/2}\leq\frac{m^{1/2}\ell^{(p+1)/2}}{\sqrt{p+1}}S_k^{1/2}.\label{eq:half-path}
\end{align}
\end{lemma}
\begin{proof}
For $t\geq1$, since $x_{k,0}=z_k$ and $x_{k,j+1}=x_{k,j}+s_{k,j}$,
\begin{equation*}
\|x_{k,t}-z_k\|\leq\sum_{j=0}^{t-1}\|s_{k,j}\|.
\end{equation*}
For nonnegative numbers $a_0,\ldots,a_{h-1}$, H\"older's inequality gives
\begin{equation}\label{eq:partial-sum-holder}
\left(\sum_{j=0}^{h-1}a_j\right)^{p+1}\leq h^p\sum_{j=0}^{h-1}a_j^{p+1}.
\end{equation}
Applying \eqref{eq:partial-sum-holder} with $h=t$ and summing over $t$ gives
\begin{equation*}
\sum_{t=0}^{m-1}\|x_{k,t}-z_k\|^{p+1}\leq\left(\sum_{a=1}^{m-1}a^p\right)S_k.
\end{equation*}
Since
\begin{equation*}
\sum_{a=1}^{h-1}a^p\leq\int_0^h u^p\,du=\frac{h^{p+1}}{p+1},
\end{equation*}
we obtain \eqref{eq:block-path}.

The same argument applies within each subblock. Indeed, if $b$ is a subblock start and $t$ lies in that subblock, then
\begin{equation*}
\|x_{k,t}-x_{k,b}\|\leq\sum_{j=b}^{t-1}\|s_{k,j}\|.
\end{equation*}
Since every subblock has length at most $\ell$, applying \eqref{eq:partial-sum-holder} within each subblock and summing over all subblocks yields
\begin{equation*}
\sum_{t=0}^{m-1}\|x_{k,t}-x_{k,b(t)}\|^{p+1}\leq\frac{\ell^{p+1}}{p+1}S_k,
\end{equation*}
which proves \eqref{eq:subblock-path}. The same argument also covers a truncated final subblock.

For \eqref{eq:mixed-path}, the Cauchy--Schwarz inequality and \eqref{eq:block-path}--\eqref{eq:subblock-path} give
\begin{align*}
\sum_{t=0}^{m-1}\left(\|x_{k,t}-z_k\|\|x_{k,t}-x_{k,b(t)}\|\right)^{(p+1)/2}
&\leq\left(\sum_{t=0}^{m-1}\|x_{k,t}-z_k\|^{p+1}\right)^{1/2}\left(\sum_{t=0}^{m-1}\|x_{k,t}-x_{k,b(t)}\|^{p+1}\right)^{1/2}\\
&\leq\frac{(m\ell)^{(p+1)/2}}{p+1}S_k.
\end{align*}

Finally, another application of Cauchy--Schwarz gives
\begin{align*}
\sum_{t=0}^{m-1}\|x_{k,t}-x_{k,b(t)}\|^{(p+1)/2}
\leq m^{1/2}\left(\sum_{t=0}^{m-1}\|x_{k,t}-x_{k,b(t)}\|^{p+1}\right)^{1/2}\leq\frac{m^{1/2}\ell^{(p+1)/2}}{\sqrt{p+1}}S_k^{1/2},
\end{align*}
which proves \eqref{eq:half-path}.
\end{proof}

\section{Technical proofs}
\label{sec:app:proofs}

\begin{proof}[\textbf{Proof of Lemma~\ref{lem:local-step-estimates}}]
Fix $t\in\{0,\ldots,m-1\}$ and write $r:=\|s_{k,t}\|$. Since $s_{k,t}$ satisfies \eqref{eq:model-decrease},
\[
f(x_{k,t})-f(x_{k,t+1})>M_{k,t}(s_{k,t})-f(x_{k,t+1}).
\]
Applying \eqref{eq:approx-value-error} with $s=s_{k,t}$ gives
\[
f(x_{k,t})-f(x_{k,t+1})\geq\frac{\sigma_k-L_p}{(p+1)!}r^{p+1}-\frac{\delta_{A,k,t}}{(p-1)!}r^{p-1}-\frac{\delta_{B,k,t}}{p!}r^p.
\]
Since \eqref{eq:sigma-large} gives $\sigma_k\geq2L_p$, we obtain
\begin{equation}\label{eq:app-decrease-before-young}
f(x_{k,t})-f(x_{k,t+1})\geq\frac{\sigma_k}{2(p+1)!}r^{p+1}-\frac{\delta_{A,k,t}}{(p-1)!}r^{p-1}-\frac{\delta_{B,k,t}}{p!}r^p.
\end{equation}
We next control the last two terms by scaled Young inequalities. For the order-$(p-1)$ term, use the conjugate exponents
$\frac{p+1}{p-1}$ and $\frac{p+1}{2}$, and for the order-$p$ term use $\frac{p+1}{p}$ and $\beta_B:=p+1$. 
Applying Young's inequality with these conjugate exponents and scaling the first terms by $\sigma_k/[8(p+1)!]$ gives
\begin{equation}\label{eq:app-young-A}
\frac{\delta_{A,k,t}}{(p-1)!}r^{p-1}\leq\frac{\sigma_k}{8(p+1)!}r^{p+1}+K_{A,p}\frac{\delta_{A,k,t}^{(p+1)/2}}{\sigma_k^{(p-1)/2}},
\end{equation}
and
\begin{equation}\label{eq:app-young-B}
\frac{\delta_{B,k,t}}{p!}r^p\leq\frac{\sigma_k}{8(p+1)!}r^{p+1}+K_{B,p}\frac{\delta_{B,k,t}^{p+1}}{\sigma_k^p}.
\end{equation}
A direct substitution in Young's inequality gives the constants $K_{A,p}$ and $K_{B,p}$ defined in Appendix~\ref{app:constants}. Substituting \eqref{eq:app-young-A}--\eqref{eq:app-young-B} into \eqref{eq:app-decrease-before-young} yields
\[
f(x_{k,t})-f(x_{k,t+1})\geq\frac{\sigma_k}{4(p+1)!}r^{p+1}-K_{A,p}\frac{\delta_{A,k,t}^{(p+1)/2}}{\sigma_k^{(p-1)/2}}-K_{B,p}\frac{\delta_{B,k,t}^{p+1}}{\sigma_k^p}.
\]
Therefore, with $c_f:=1/[4(p+1)!]$, $K_{f,p}:=\max\{K_{A,p},K_{B,p}\}$, and $\mathcal E_{k,t}$ defined in \eqref{eq:local-error-penalty},
\[
f(x_{k,t})-f(x_{k,t+1})\geq c_f\sigma_k\|s_{k,t}\|^{p+1}-K_{f,p}\mathcal E_{k,t},
\]
which proves \eqref{eq:local-decrease}.
We now prove \eqref{eq:local-gradient}. The residual condition \eqref{eq:model-residual} gives
\[
\|\nabla_sM_{k,t}(s_{k,t})\|\leq\theta\frac{\sigma_k}{p!}r^p.
\]
Combining this inequality with \eqref{eq:approx-gradient-error} gives
\[
\|\nabla f(x_{k,t+1})\|\leq\frac{\theta\sigma_k+\sigma_k+L_p}{p!}r^p+\frac{\delta_{A,k,t}}{(p-2)!}r^{p-2}+\frac{\delta_{B,k,t}}{(p-1)!}r^{p-1}.
\]
Since $L_p\leq\sigma_k$, it follows that
\begin{equation}\label{eq:app-gradient-before-power}
\|\nabla f(x_{k,t+1})\|\leq a_0\sigma_kr^p+a_A\delta_{A,k,t}r^{p-2}+a_B\delta_{B,k,t}r^{p-1},
\end{equation}
where
\[
a_0:=\frac{\theta+2}{p!},\qquad a_A:=\frac{1}{(p-2)!},\qquad a_B:=\frac{1}{(p-1)!}.
\]
Set $\tau:=(p+1)/p$. Since $\tau\geq1$,
\begin{equation}\label{eq:app-three-term-power}
(u+v+w)^\tau\leq3^{\tau-1}(u^\tau+v^\tau+w^\tau)=3^{1/p}(u^\tau+v^\tau+w^\tau)
\end{equation}
for all $u,v,w\geq0$. Moreover,
\begin{equation}\label{eq:app-gradient-energy}
(\sigma_kr^p)^\tau=\sigma_k^{1/p}\left(\sigma_kr^{p+1}\right).
\end{equation}
For the order-$(p-1)$ error term, weighted Young's inequality gives
\begin{equation}\label{eq:app-gradient-young-A}
\left(\delta_{A,k,t}r^{p-2}\right)^\tau\leq\sigma_k^{1/p}\left[\frac{p-2}{p}\sigma_kr^{p+1}+\frac{2}{p}\frac{\delta_{A,k,t}^{(p+1)/2}}{\sigma_k^{(p-1)/2}}\right].
\end{equation}
Similarly,
\begin{equation}\label{eq:app-gradient-young-B}
\left(\delta_{B,k,t}r^{p-1}\right)^\tau\leq\sigma_k^{1/p}\left[\frac{p-1}{p}\sigma_kr^{p+1}+\frac{1}{p}\frac{\delta_{B,k,t}^{p+1}}{\sigma_k^p}\right].
\end{equation}
Raising \eqref{eq:app-gradient-before-power} to the power $\tau$, applying \eqref{eq:app-three-term-power}, and then using \eqref{eq:app-gradient-energy}--\eqref{eq:app-gradient-young-B} yields
\begin{align*}
\|\nabla f(x_{k,t+1})\|^\tau\leq3^{1/p}\sigma_k^{1/p}\Bigg[\left(a_0^\tau+\frac{p-2}{p}a_A^\tau+\frac{p-1}{p}a_B^\tau\right)\sigma_kr^{p+1}
+\frac{2}{p}a_A^\tau\frac{\delta_{A,k,t}^{(p+1)/2}}{\sigma_k^{(p-1)/2}}+\frac{1}{p}a_B^\tau\frac{\delta_{B,k,t}^{p+1}}{\sigma_k^p}\Bigg].
\end{align*}
By the definition of $C_{g,p}$ in Appendix~\ref{app:constants}, the three coefficients in brackets are bounded by $C_{g,p}/3^{1/p}$. Hence
\[
\|\nabla f(x_{k,t+1})\|^{(p+1)/p}\leq C_{g,p}\sigma_k^{1/p}\left(\sigma_k\|s_{k,t}\|^{p+1}+\mathcal E_{k,t}\right),
\]
which proves \eqref{eq:local-gradient}.
\end{proof}

\begin{proof}[\textbf{Proof of Lemma~\ref{lem:accumulated-error}}]
For brevity, set
\[
D_t:=\|x_{k,t}-z_k\|,\qquad E_t:=\|x_{k,t}-x_{k,b(t)}\|.
\]
Since $2R=p+1$, the bounds \eqref{eq:B-error-along-block}--\eqref{eq:A-error-along-block} give
\begin{align}\label{eq:app-A-error-sum}
\sum_{t=0}^{m-1}\frac{\delta_{A,k,t}^{R}}{\sigma_k^{R-1}}\leq H_p\Bigg[\frac{m\ov\delta_{A,k}^{R}}{\sigma_k^{R-1}}+\frac{\ov\delta_{B,k}^{R}}{\sigma_k^{R-1}}\sum_{t=0}^{m-1}E_t^{R}
+\frac{L_p^{R}}{\sigma_k^{R-1}}\sum_{t=0}^{m-1}(D_tE_t)^{R}+\frac{(L_p/2)^{R}}{\sigma_k^{R-1}}\sum_{t=0}^{m-1}E_t^{p+1}\Bigg],
\end{align}
where we used $(a_1+a_2+a_3+a_4)^R\leq4^{R-1}\sum_{i=1}^4a_i^R$. Similarly,
\begin{equation}\label{eq:app-B-error-sum}
\sum_{t=0}^{m-1}\frac{\delta_{B,k,t}^{p+1}}{\sigma_k^p}\leq \frac{2^pm\ov\delta_{B,k}^{p+1}}{\sigma_k^p}+\frac{2^pL_p^{p+1}}{\sigma_k^p}\sum_{t=0}^{m-1}D_t^{p+1}.
\end{equation}
We first control the path-dependent terms. By \eqref{eq:mixed-path},
\begin{align*}
\frac{H_pL_p^R}{\sigma_k^{R-1}}\sum_{t=0}^{m-1}(D_tE_t)^R
\leq\frac{H_pL_p^R(m\ell)^R}{(p+1)\sigma_k^{R-1}}S_k
=\frac{H_p}{p+1}\left(\frac{L_pm\ell}{\sigma_k}\right)^R\sigma_kS_k\leq\frac{H_p}{(p+1)C_\sigma^R}\sigma_kS_k
\leq\frac{\rho}{4}\sigma_kS_k,
\end{align*}
where we used $L_p\leq L_k$, $\sigma_k=C_\sigma L_km\ell$, and the choice of $C_\sigma$ in Appendix~\ref{app:constants}. 
Next, \eqref{eq:subblock-path} gives
\begin{align*}
\frac{H_p(L_p/2)^R}{\sigma_k^{R-1}}\sum_{t=0}^{m-1}E_t^{p+1}
\leq\frac{H_p(L_p/2)^R\ell^{p+1}}{(p+1)\sigma_k^{R-1}}S_k
=\frac{H_p}{2^R(p+1)}\left(\frac{L_p\ell^2}{\sigma_k}\right)^R\sigma_kS_k
\leq\frac{H_p}{2^R(p+1)C_\sigma^R}\left(\frac{\ell}{m}\right)^R\sigma_kS_k
\leq\frac{\rho}{4}\sigma_kS_k.
\end{align*}
Here we also used $\ell\leq m$.
For the path-dependent part of the order-$p$ error, \eqref{eq:block-path} yields
\begin{align*}
\frac{2^pL_p^{p+1}}{\sigma_k^p}\sum_{t=0}^{m-1}D_t^{p+1}
\leq\frac{2^pL_p^{p+1}m^{p+1}}{(p+1)\sigma_k^p}S_k
=\frac{2^p}{p+1}\left(\frac{L_pm}{\sigma_k}\right)^{p+1}\sigma_kS_k
\leq\frac{2^p}{(p+1)C_\sigma^{p+1}}\sigma_kS_k
\leq\frac{\rho}{4}\sigma_kS_k.
\end{align*}
It remains to control the term involving the fresh order-$p$ reconstruction error and the subblock displacement. By \eqref{eq:half-path},
\begin{align*}
\frac{H_p\ov\delta_{B,k}^{R}}{\sigma_k^{R-1}}\sum_{t=0}^{m-1}E_t^R
&\leq\frac{H_p\ov\delta_{B,k}^{R}m^{1/2}\ell^R}{\sqrt{p+1}\,\sigma_k^{R-1}}S_k^{1/2}.
\end{align*}
Writing the right-hand side as $uv$ with
\[
u:=(\sigma_kS_k)^{1/2},\qquad v:=\frac{H_pm^{1/2}\ell^R\ov\delta_{B,k}^{R}}{\sqrt{p+1}\,\sigma_k^{R-1/2}},
\]
and applying $uv\leq(\rho/4)u^2+v^2/\rho$, we obtain
\begin{equation}
\frac{H_p\ov\delta_{B,k}^{R}}{\sigma_k^{R-1}}\sum_{t=0}^{m-1}E_t^R
\leq\frac{\rho}{4}\sigma_kS_k+\frac{H_p^2}{(p+1)\rho}\frac{m\ell^{p+1}\ov\delta_{B,k}^{p+1}}{\sigma_k^p}.
\label{eq:app-B-transport-control}
\end{equation}
Combining \eqref{eq:app-A-error-sum}--\eqref{eq:app-B-transport-control}, the four path-dependent contributions sum to at most $\rho\sigma_kS_k$. Therefore,
\begin{align}
\sum_{t=0}^{m-1}\mathcal E_{k,t}
\leq{}&
\rho\sigma_kS_k
+\frac{H_pm\ov\delta_{A,k}^{R}}{\sigma_k^{R-1}}
+\frac{2^pm\ov\delta_{B,k}^{p+1}}{\sigma_k^p}
+\frac{H_p^2}{(p+1)\rho}\frac{m\ell^{p+1}\ov\delta_{B,k}^{p+1}}{\sigma_k^p}.
\label{eq:app-error-before-fresh}
\end{align}
It remains to bound the three terms involving the fresh reconstruction errors. By \eqref{eq:delta-bars},
\begin{equation}
\frac{\ov\delta_{A,k}^{R}}{\sigma_k^{R-1}}
=
c_A^R\eps^{(p+1)/p}\sigma_k^{-1/p},
\label{eq:app-fresh-A}
\end{equation}
and
\begin{equation}
\frac{\ov\delta_{B,k}^{p+1}}{\sigma_k^p}
=
c_B^{p+1}\ell^{-(p+1)}\eps^{(p+1)/p}\sigma_k^{-1/p}.
\label{eq:app-fresh-B}
\end{equation}
Substituting \eqref{eq:app-fresh-A}--\eqref{eq:app-fresh-B} into \eqref{eq:app-error-before-fresh} and using $\ell\geq1$ gives
\begin{align*}
\sum_{t=0}^{m-1}\mathcal E_{k,t}\leq\rho\sigma_kS_k+m\eps^{(p+1)/p}\sigma_k^{-1/p}\left[H_pc_A^R+\left(2^p+\frac{H_p^2}{(p+1)\rho}\right)c_B^{p+1}\right].
\end{align*}
The choices of $c_A$ and $c_B$ in Appendix~\ref{app:constants} satisfy
\[
H_pc_A^R\leq\frac{\eta_{\rm err}}{2},\qquad\left(2^p+\frac{H_p^2}{(p+1)\rho}\right)c_B^{p+1}\leq\frac{\eta_{\rm err}}{2}.
\]
Consequently,
\[
\sum_{t=0}^{m-1}\mathcal E_{k,t}\leq\rho\sigma_k\sum_{t=0}^{m-1}\|s_{k,t}\|^{p+1}+\eta_{\rm err}m\eps^{(p+1)/p}\sigma_k^{-1/p},
\]
which proves \eqref{eq:accumulated-error}.
\end{proof}

\begin{proof}[\textbf{Proof of Proposition~\ref{prop:finite-inner-realization}}]
Since $\nabla_sM(0)\neq0$, the directional derivative of $M$ at the origin along $-\nabla_sM(0)$ is strictly negative. Hence the first backtracking search accepts a positive stepsize and produces
$M(s_1)<M(0)$. All subsequent accepted iterates belong to the sublevel set
\[
\mathcal L_1 := \{s\in\R^n:M(s)\leq M(s_1)\}.
\]
Since $M$ is coercive, $\mathcal L_1$ is compact. Moreover, $M(s_1)<M(0)$ implies $0\notin\mathcal L_1$. Hence there exists $d_1>0$ such that
\begin{equation}\label{eq:inner-away-zero}
\|s\|\geq d_1, \qquad \forall s\in\mathcal L_1.
\end{equation}
Consider
\[
\mathcal K :=\left\{s-\alpha\nabla_sM(s):s\in\mathcal L_1,\; 0\leq\alpha\leq\alpha_{\max}
\right\}.
\]
The set $\mathcal K$ is compact. Since $M$ is twice continuously differentiable, $\nabla_sM$ is Lipschitz continuous on the compact set $\operatorname{conv}(\mathcal K)$ with some constant $L_M>0$. Therefore, for every $s\in\mathcal L_1$ and
$0\leq\alpha\leq\alpha_{\max}$,
\[
M(s-\alpha\nabla_sM(s)) \leq M(s)-\alpha\left(1-\frac{L_M\alpha}{2}\right)\|\nabla_sM(s)\|^2.
\]
Thus \eqref{eq:inner-armijo} holds whenever 
$\alpha \leq \frac{2(1-c)}{L_M}$.
By the backtracking rule,
\begin{equation}\label{eq:inner-alpha-lower}
\alpha_j\geq\beta\min\left\{\alpha_{\max},\frac{2(1-c)}{L_M}\right\}=:\alpha_{\min}>0.
\end{equation}
Using \eqref{eq:inner-armijo},
\[
M(s_{j+1}) \leq M(s_j)-c\alpha_j \|\nabla_s M(s_j)\|^2.
\]
Since $\{M(s_j)\}$ is bounded below and $\alpha_j\geq\alpha_{\min}$, summation gives
\[
\sum_{j=1}^{\infty}\|\nabla_sM(s_j)\|^2 <\infty.
\]
Hence $\|\nabla_sM(s_j)\|\to 0$. On the other hand, \eqref{eq:inner-away-zero} gives
\[
\theta\frac{\sigma}{p!}\|s_j\|^p \geq \theta\frac{\sigma}{p!}d_1^p>0,
\]
for every $j\geq 1$. Therefore, after finitely many accepted steps,
\[
\|\nabla_sM(s_j)\|\leq \theta\frac{\sigma}{p!}\|s_j\|^p.
\]
Moreover, the Armijo iteration preserves
\[
M(s_j)\leq M(s_1)<M(0),
\]
for every $j\geq 1$. Thus both conditions in \eqref{eq:inner-admissible} hold after finitely
many accepted steps.
\end{proof}

\bibliographystyle{spbasic}
\bibliography{references}

\begin{thebibliography}{23}
\providecommand{\natexlab}[1]{#1}
\providecommand{\url}[1]{{#1}}
\providecommand{\urlprefix}{URL }
\expandafter\ifx\csname urlstyle\endcsname\relax
  \providecommand{\doi}[1]{DOI~\discretionary{}{}{}#1}\else
  \providecommand{\doi}{DOI~\discretionary{}{}{}\begingroup
  \urlstyle{rm}\Url}\fi
\providecommand{\eprint}[2][]{\url{#2}}

\bibitem[{Bellavia et~al(2019)Bellavia, Gurioli, Morini, and
  Toint}]{BellaviaEtAl2019}
Bellavia S, Gurioli G, Morini B, Toint PL (2019) Adaptive regularization
  algorithms with inexact evaluations for nonconvex optimization. SIAM Journal
  on Optimization 29(4):2881--2915

\bibitem[{Birgin et~al(2017)Birgin, Gardenghi, Mart\'inez, Santos, and
  Toint}]{BirginEtAl2017}
Birgin EG, Gardenghi JL, Mart\'inez JM, Santos SA, Toint PL (2017) Worst-case
  evaluation complexity for unconstrained nonlinear optimization using
  high-order regularized models. Mathematical Programming 163:359--368

\bibitem[{Carmon et~al(2020)Carmon, Duchi, Hinder, and
  Sidford}]{CarmonEtAl2020}
Carmon Y, Duchi JC, Hinder O, Sidford A (2020) Lower bounds for finding
  stationary points {I}. Mathematical Programming 184:71--120

\bibitem[{Cartis and Zhu(2026{\natexlab{a}})}]{CartisZhu2025}
Cartis C, Zhu W (2026{\natexlab{a}}) Cubic-quartic regularization models for
  solving polynomial subproblems in third-order tensor methods. Mathematical
  Programming 217:167--219

\bibitem[{Cartis and Zhu(2026{\natexlab{b}})}]{CartisZhu2026}
Cartis C, Zhu W (2026{\natexlab{b}}) Second-order methods for
  quartically-regularised cubic polynomials, with applications to high-order
  tensor methods. Mathematical Programming 215:669--715

\bibitem[{Cartis et~al(2011)Cartis, Gould, and Toint}]{CartisGouldToint2011ARC}
Cartis C, Gould NIM, Toint PL (2011) Adaptive cubic regularisation methods for
  unconstrained optimization. part {II}: Worst-case function- and
  derivative-evaluation complexity. Mathematical Programming 130:295--319

\bibitem[{Cartis et~al(2020)Cartis, Gould, and
  Toint}]{CartisGouldToint2020Sharp}
Cartis C, Gould NIM, Toint PL (2020) Sharp worst-case evaluation complexity
  bounds for arbitrary-order nonconvex optimization with inexpensive
  constraints. SIAM Journal on Optimization 30(1):513--541

\bibitem[{Cartis et~al(2025)Cartis, Jerad, and Welzel}]{CartisJerad2025}
Cartis C, Jerad S, Welzel K (2025) On global rates for regularization methods
  based on secant derivative approximations. \doi{10.48550/arXiv.2509.07580},
  arXiv:2509.07580, \eprint{2509.07580}

\bibitem[{Cartis et~al(2026)Cartis, Hauser, Liu, Welzel, and
  Zhu}]{CartisEtAl2026AR3}
Cartis C, Hauser R, Liu Y, Welzel K, Zhu W (2026) Efficient implementation of
  third-order tensor methods with adaptive regularization for unconstrained
  optimization. Mathematical Programming Computation 18:877--948

\bibitem[{Doikov and Grapiglia(2025)}]{DoikovGrapiglia2025}
Doikov N, Grapiglia GN (2025) First and zeroth-order implementations of the
  regularized {N}ewton method with lazy approximated {H}essians. Journal of
  Scientific Computing 103:32

\bibitem[{Doikov and Grapiglia(2026)}]{DoikovGrapiglia2026}
Doikov N, Grapiglia GN (2026) On the complexity of lower-order implementations
  of higher-order methods. Mathematical Programming

\bibitem[{Doikov et~al(2023)Doikov, Chayti, and Jaggi}]{DoikovChaytiJaggi2023}
Doikov N, Chayti EM, Jaggi M (2023) Second-order optimization with lazy
  {H}essians. In: Proceedings of the 40th International Conference on Machine
  Learning, PMLR, Proceedings of Machine Learning Research, vol 202, pp
  8138--8161

\bibitem[{Gon{\c{c}}alves and Filho(2026)}]{GoncalvesGehlen2026}
Gon{\c{c}}alves MLN, Filho VG (2026) Inexact cubic regularization method with
  adaptive reuse of {H}essian approximations. Optimization Online

\bibitem[{Grapiglia(2023)}]{Grapiglia2023}
Grapiglia GN (2023) Quadratic regularization methods with finite-difference
  gradient approximations. Computational Optimization and Applications
  85:683--703

\bibitem[{Grapiglia et~al(2022)Grapiglia, Gon\c{c}alves, and
  Silva}]{GrapigliaGoncalvesSilva2022}
Grapiglia GN, Gon\c{c}alves MLN, Silva GN (2022) A cubic regularization of
  {N}ewton's method with finite difference {H}essian approximations. Numerical
  Algorithms 90:607--630

\bibitem[{Griewank(1981)}]{griewank1981modification}
Griewank A (1981) The modification of {N}ewton's method for unconstrained
  optimization by bounding cubic terms. Tech. Rep. NA/12, Department of Applied
  Mathematics and Theoretical Physics, University of Cambridge

\bibitem[{Kabgani and Ahookhosh(2026)}]{Kabgani2026Tensor}
Kabgani A, Ahookhosh M (2026) Tensor envelopes in the nonconvex setting.
  \doi{10.48550/arXiv.2610.04610}, \eprint{2610.04610}

\bibitem[{Nesterov(2018)}]{Nesterov2018}
Nesterov Y (2018) Lectures on Convex Optimization, 2nd edn. Springer Cham

\bibitem[{Nesterov(2021)}]{Nesterov2021Superfast}
Nesterov Y (2021) Superfast second-order methods for unconstrained convex
  optimization. Journal of Optimization Theory and Applications 191:1--30

\bibitem[{Nesterov and Polyak(2006)}]{NesterovPolyak2006}
Nesterov Y, Polyak BT (2006) Cubic regularization of {N}ewton method and its
  global performance. Mathematical Programming 108:177--205

\bibitem[{Thomas(2014)}]{Thomas2014Polarization}
Thomas EGF (2014) A polarization identity for multilinear maps. Indagationes
  Mathematicae 25(3):468--474

\bibitem[{Welzel and Hauser(2024)}]{WelzelHauser2024}
Welzel K, Hauser RA (2024) Approximating higher-order derivative tensors using
  secant updates. SIAM Journal on Optimization 34(1):893--917

\bibitem[{Zhou(2026)}]{Zhou2026FirstOrder}
Zhou D (2026) Sharp first-order lower bounds for higher-order smooth nonconvex
  optimization. \doi{10.48550/arXiv.2606.05438}, \eprint{2606.05438}

\end{thebibliography}

\end{document}